\documentclass[a4paper]{article}

\usepackage{MyLaTeX,MyArticle}
\title{The structure of projective indecomposable modules for $A_n$, $n\leq 12$}
\author{David A.~Craven, University of Birmingham}
\date{May 2016}

\begin{document}
\maketitle

\begin{abstract}
This article gives the structure of all projective indecomposable modules in blocks with non-cyclic defect group for $n\leq 11$, and almost all in the case of $n=12$, leaving four simple modules for $p=2$ and one simple module for $p=3$.
\end{abstract}

In the representation theory of the symmetric and alternating groups, blocks with abelian defect group are, in some sense, well understood. The decomposition numbers are not known in general, the $\Ext^1$ matrices are not known in general, the structure of the projectives are not known in general, but things are better than when the defect group is non-abelian. A lack of explicit examples is of no help, of course, and the purpose of this short paper is to provide researchers with a few more examples of $\Ext^1$-matrices and the socle structures of the projective indecomposable modules for $A_n$ for $n\leq 12$. This bound is largely imposed on the author by the fact that $A_{12}$ is already difficult to do via computer without advanced methods. Indeed, it is not the construction of the PIMs but analysing their socle structure that is the real issue, taking many times longer than their construction.

Although we are mainly interested in non-abelian defect groups -- so $p=2,3$ -- we include the $\Ext^1$-matrices and PIM structures for $p=5$ as well, since they might be of use to people. For $p\geq 7$ the defect group is cyclic and the structure of all blocks are known by the theory of Brauer trees, as is the case for $p=5$ and $n\leq 9$. We omit the description of all projectives that lie in blocks of defect $0$ or $1$.

If we look through the literature, we can find the following cases where the projective indecomposable modules are already described. There might well be others, but I cannot find them at the moment.

\begin{center}
\begin{tabular}{cccc}
\hline $n$ & $A_n$, $S_n$, $p=2$ & $A_n$, $p=3$ & $S_n$, $p=3$
\\\hline $n=7$ & Erdmann \cite{erdmann1977} & Scopes \cite{scopes1995} & Scopes \cite{scopes1995}
\\ $n=8$ & Benson \cite{benson1983} & Scopes \cite{scopes1995} & Scopes \cite{scopes1995} 
\\ $n=9$ & Benson \cite{benson1983a} & Siegel \cite{siegel1991} & Tan \cite{tan2000}
\\ $n=10$ & All here & Tan \cite{tan2000} & Tan \cite{tan2000}
\\ $n=11$ & All here & Three modules by Tan \cite{tan2000a}, rest here & Six modules by Tan \cite{tan2000a}, rest here
\\ $n=12$ & All but four here & All but one here & This paper
\\ \hline
\end{tabular}
\end{center}

Notice that for $p=2$ we simply induce from $A_n$ to $S_n$, replacing a module $M$ in the diagram for each projective by the self-extension of $M$. (This also allows us to go backwards, of course.) For $S_{10}$ and $p=3$, the projectives are described in \cite{tan2000}, and those for $A_{10}$ can be obtained from them by restriction, although we reproduce them here for ease of reference. For $p=5$ the blocks have weight $2$, and are completely understood by work of Scopes \cite{scopes1995}. We give the projectives here just to make life easier for people who want them.

\section{$A_{10}$ and $S_{10}$, characteristic $2$}

We firstly compute the projective indecomposable modules for the alternating group $A_{10}$, and then induce them to $S_{10}$ to produce the projectives there.

There are seven simple modules in the principal block of $S_{10}$ (and $A_{10}$) and they are below.
\begin{center}
\begin{tabular}{ccc}
\hline $\lambda$ &$\dim(D^\lambda)$ & Factors of $S^\lambda$
\\\hline  $(10)$ & $1$ & $1$
\\ $(9,1)$ & $8$ & $8,1$
\\ $(6,4)$ & $16$ & $16,48,26$
\\ $(8,2)$ & $26$ & $26,8,1$
\\ $(7,3)$ & $48$ & $48,26,1$
\\ $(6,3,1)$ & $198$ & $198,16,48,1,26^2$
\\ $(5,3,2)$ & $200$ & $200,198,16,26,8,1^2$
\\ \hline
\end{tabular}
\end{center}

The non-principal block that does not have defect zero has two characters, $S^{(7,3,1)}=D^{(7,2,1)}$ of degree $160$ and $D^{(5,4,1)}$ of degree $128$.

Here are the projectives for the alternating group $A_{10}$, which we can then induce to $S_{10}$. For completeness we list them as well.

\[\begin{array}{c}
1
\\8\;\;26\;\;200
\\1\;\;1\;\;1\;\;48\;\;198
\\8\;\;8\;\;26\;\;26\;\;26\;\;200
\\1\;\;1\;\;1\;\;1\;\;1\;\;1\;\;16\;\;16\;\;16\;\;48
\\8\;\;8\;\;8\;\;26\;\;26\;\;26\;\;26\;\;26\;\;198\;\;200
\\1\;\;1\;\;1\;\;1\;\;1\;\;1\;\;1\;\;16\;\;16\;\;48\;\;48\;\;48\;\;198
\\8\;\;8\;\;8\;\;26\;\;26\;\;26\;\;26\;\;26\;\;26\;\;26\;\;198\;\;200\;\;200\;\;200
\\1\;\;1\;\;1\;\;1\;\;1\;\;1\;\;1\;\;1\;\;1\;\;16\;\;16\;\;16\;\;48\;\;48\;\;198
\\1\;\;8\;\;8\;\;8\;\;8\;\;8\;\;26\;\;26\;\;26\;\;26\;\;26\;\;26\;\;48\;\;198\;\;200\;\;200
\\1\;\;1\;\;1\;\;1\;\;1\;\;1\;\;1\;\;1\;\;1\;\;8\;\;16\;\;16\;\;16\;\;16\;\;16\;\;26\;\;48\;\;48\;\;198
\\1\;\;1\;\;8\;\;8\;\;8\;\;16\;\;26\;\;26\;\;26\;\;26\;\;26\;\;26\;\;26\;\;200\;\;200
\\1\;\;1\;\;1\;\;1\;\;1\;\;1\;\;16\;\;16\;\;26\;\;48\;\;48\;\;48\;\;198
\\8\;\;8\;\;8\;\;26\;\;26\;\;26\;\;26\;\;48\;\;198\;\;200
\\1\;\;1\;\;1\;\;1\;\;1\;\;16\;\;16\;\;26\;\;48
\\8\;\;8\;\;16\;\;26\;\;26\;\;198\;\;200\;\;200
\\1\;\;1\;\;1\;\;1\;\;1\;\;1\;\;16\;\;26\;\;48
\\8\;\;26\;\;198\;\;198\;\;200
\\1
\end{array} %\]
\quad
%\[
\begin{array}{c}
8
\\1
\\26\;\;200
\\1\;\;1\;\;48
\\8\;\;8\;\;26
\\1\;\;1\;\;1\;\;16\;\;16
\\8\;\;26\;\;26\;\;26\;\;200
\\1\;\;1\;\;1\;\;48\;\;48\;\;198
\\8\;\;8\;\;26\;\;26\;\;26\;\;198\;\;200
\\1\;\;1\;\;1\;\;1\;\;1\;\;16\;\;16\;\;16
\\1\;\;8\;\;8\;\;26\;\;26\;\;26\;\;200
\\1\;\;1\;\;1\;\;16\;\;48\;\;48\;\;198
\\8\;\;26\;\;26\;\;26\;\;200
\\1\;\;1\;\;1\;\;16
\\8\;\;8\;\;26\;\;198
\\1\;\;1\;\;16\;\;48
\\26\;\;200
\\1
\\8
\end{array}\]
\[\begin{array}{c}
16
\\26
\\48\;\;198
\\26\;\;26\;\;200
\\1\;\;1\;\;1\;\;16
\\8\;\;8\;\;26
\\1\;\;1\;\;16\;\;16\;\;48
\\26\;\;26\;\;26\;\;200\;\;200
\\1\;\;1\;\;1\;\;48\;\;48\;\;198
\\8\;\;8\;\;8\;\;26\;\;26\;\;26\;\;198
\\1\;\;1\;\;1\;\;1\;\;1\;\;16\;\;16\;\;16\;\;16
\\1\;\;8\;\;26\;\;26\;\;26\;\;26\;\;200
\\1\;\;1\;\;48\;\;48\;\;48\;\;198
\\8\;\;26\;\;26\;\;26
\\1\;\;1\;\;16\;\;16\;\;16\;\;198
\\1\;\;8\;\;26\;\;26\;\;200
\\1\;\;16\;\;48\;\;198
\\26\;\;200
\\16
\end{array} %\]
\quad
%\[
\begin{array}{c}
26
\\1\;\;16\;\;48\;\;198
\\8\;\;26\;\;26\;\;26
\\1\;\;1\;\;1\;\;16\;\;16\;\;48\;\;198
\\8\;\;26\;\;26\;\;26\;\;26\;\;26\;\;200
\\1\;\;1\;\;1\;\;1\;\;1\;\;16\;\;48\;\;48\;\;198
\\8\;\;8\;\;8\;\;26\;\;26\;\;26\;\;26\;\;200
\\1\;\;1\;\;1\;\;1\;\;1\;\;1\;\;1\;\;16\;\;16\;\;16\;\;48\;\;48
\\8\;\;8\;\;8\;\;26\;\;26\;\;26\;\;26\;\;26\;\;198\;\;200\;\;200
\\1\;\;1\;\;1\;\;1\;\;1\;\;1\;\;16\;\;16\;\;16\;\;48\;\;48\;\;48\;\;198
\\1\;\;8\;\;8\;\;8\;\;26\;\;26\;\;26\;\;26\;\;26\;\;26\;\;26\;\;198\;\;200\;\;200
\\1\;\;1\;\;1\;\;1\;\;1\;\;1\;\;1\;\;16\;\;16\;\;16\;\;16\;\;48\;\;48\;\;198
\\1\;\;8\;\;8\;\;8\;\;26\;\;26\;\;26\;\;26\;\;26\;\;26\;\;198\;\;200
\\1\;\;1\;\;1\;\;1\;\;16\;\;16\;\;16\;\;48\;\;48\;\;48\;\;198
\\1\;\;8\;\;26\;\;26\;\;26\;\;26\;\;26\;\;200
\\1\;\;1\;\;16\;\;16\;\;48\;\;198\;\;198
\\1\;\;8\;\;26\;\;26\;\;26
\\1\;\;16\;\;48\;\;198
\\26
\end{array}\]

\[\begin{array}{c}
48
\\26
\\1\;\;16
\\8\;\;26
\\1\;\;48\;\;198
\\26\;\;26\;\;200
\\1\;\;1\;\;1\;\;16\;\;48
\\8\;\;8\;\;26\;\;26
\\1\;\;1\;\;16\;\;16\;\;48
\\1\;\;26\;\;26\;\;26\;\;200
\\1\;\;1\;\;48\;\;48\;\;198
\\8\;\;8\;\;26\;\;26\;\;198
\\1\;\;1\;\;1\;\;16\;\;16\;\;16
\\1\;\;26\;\;26\;\;26\;\;200
\\1\;\;48\;\;198
\\8\;\;26
\\1\;\;16
\\26\;\;48
\\48
\end{array}
\quad
\begin{array}{c}
198
\\26
\\1\;\;16
\\26
\\48
\\1\;\;26
\\1\;\;200
\\1\;\;8\;\;198
\\1\;\;8\;\;16\;\;26
\\1\;\;16\;\;26\;\;200
\\1\;\;26\;\;48\;\;198
\\8\;\;26\;\;48
\\1\;\;16\;\;26
\\1\;\;26
\\8\;\;16\;\;48\;\;198
\\1\;\;26\;\;26
\\16\;\;198\;\;200
\\1\;\;1\;\;26
\\198
\end{array}
\quad
\begin{array}{c}
200
\\1
\\8
\\1\;\;16
\\26\;\;200
\\1\;\;48
\\8\;\;26\;\;198
\\1\;\;1\;\;1\;\;16\;\;16
\\8\;\;26\;\;26\;\;200
\\1\;\;1\;\;48\;\;198
\\8\;\;26\;\;26
\\1\;\;1\;\;16
\\8\;\;26
\\1\;\;48
\\26\;\;200
\\1\;\;1\;\;16
\\8\;\;198\;\;200
\\1\;\;16
\\200
\end{array}\]
Here are the symmetric group ones.

\[\begin{array}{c}
1
\\1\;\;8\;\;26\;\;200
\\1\;\;1\;\;1\;\;8\;\;26\;\;48\;\;198\;\;200
\\1\;\;1\;\;1\;\;8\;\;8\;\;26\;\;26\;\;26\;\;48\;\;198\;\;200
\\1\;\;1\;\;1\;\;1\;\;1\;\;1\;\;8\;\;8\;\;16\;\;16\;\;16\;\;26\;\;26\;\;26\;\;48\;\;200
\\1\;\;1\;\;1\;\;1\;\;1\;\;1\;\;8\;\;8\;\;8\;\;16\;\;16\;\;16\;\;26\;\;26\;\;26\;\;26\;\;26\;\;48\;\;198\;\;200
\\1\;\;1\;\;1\;\;1\;\;1\;\;1\;\;1\;\;8\;\;8\;\;8\;\;16\;\;16\;\;26\;\;26\;\;26\;\;26\;\;26\;\;48\;\;48\;\;48\;\;198\;\;198\;\;200
\\1\;\;1\;\;1\;\;1\;\;1\;\;1\;\;1\;\;8\;\;8\;\;8\;\;16\;\;16\;\;26\;\;26\;\;26\;\;26\;\;26\;\;26\;\;26\;\;48\;\;48\;\;48\;\;198\;\;198\;\;200\;\;200\;\;200
\\1\;\;1\;\;1\;\;1\;\;1\;\;1\;\;1\;\;1\;\;1\;\;8\;\;8\;\;8\;\;16\;\;16\;\;16\;\;26\;\;26\;\;26\;\;26\;\;26\;\;26\;\;26\;\;48\;\;48\;\;198\;\;198\;\;200\;\;200\;\;200
\\1\;\;1\;\;1\;\;1\;\;1\;\;1\;\;1\;\;1\;\;1\;\;1\;\;8\;\;8\;\;8\;\;8\;\;8\;\;16\;\;16\;\;16\;\;26\;\;26\;\;26\;\;26\;\;26\;\;26\;\;48\;\;48\;\;48\;\;198\;\;198\;\;200\;\;200
\\1\;\;1\;\;1\;\;1\;\;1\;\;1\;\;1\;\;1\;\;1\;\;1\;\;8\;\;8\;\;8\;\;8\;\;8\;\;8\;\;16\;\;16\;\;16\;\;16\;\;16\;\;26\;\;26\;\;26\;\;26\;\;26\;\;26\;\;26\;\;48\;\;48\;\;48\;\;198\;\;198\;\;200\;\;200
\\1\;\;1\;\;1\;\;1\;\;1\;\;1\;\;1\;\;1\;\;1\;\;1\;\;1\;\;8\;\;8\;\;8\;\;8\;\;16\;\;16\;\;16\;\;16\;\;16\;\;16\;\;26\;\;26\;\;26\;\;26\;\;26\;\;26\;\;26\;\;26\;\;48\;\;48\;\;198\;\;200\;\;200
\\1\;\;1\;\;1\;\;1\;\;1\;\;1\;\;1\;\;1\;\;8\;\;8\;\;8\;\;16\;\;16\;\;16\;\;26\;\;26\;\;26\;\;26\;\;26\;\;26\;\;26\;\;26\;\;48\;\;48\;\;48\;\;198\;\;200\;\;200
\\1\;\;1\;\;1\;\;1\;\;1\;\;1\;\;8\;\;8\;\;8\;\;16\;\;16\;\;26\;\;26\;\;26\;\;26\;\;26\;\;48\;\;48\;\;48\;\;48\;\;198\;\;198\;\;200
\\1\;\;1\;\;1\;\;1\;\;1\;\;8\;\;8\;\;8\;\;16\;\;16\;\;26\;\;26\;\;26\;\;26\;\;26\;\;48\;\;48\;\;198\;\;200
\\1\;\;1\;\;1\;\;1\;\;1\;\;8\;\;8\;\;16\;\;16\;\;16\;\;26\;\;26\;\;26\;\;48\;\;198\;\;200\;\;200
\\1\;\;1\;\;1\;\;1\;\;1\;\;1\;\;8\;\;8\;\;16\;\;16\;\;26\;\;26\;\;26\;\;48\;\;198\;\;200\;\;200
\\1\;\;1\;\;1\;\;1\;\;1\;\;1\;\;8\;\;16\;\;26\;\;26\;\;48\;\;198\;\;198\;\;200
\\1\;\;8\;\;26\;\;198\;\;198\;\;200
\\1
\end{array}\]

\[\begin{array}{c} 8
\\1\;\;8
\\1\;\;26\;\;200
\\1\;\;1\;\;26\;\;48\;\;200
\\1\;\;1\;\;8\;\;8\;\;26\;\;48
\\1\;\;1\;\;1\;\;8\;\;8\;\;16\;\;16\;\;26
\\1\;\;1\;\;1\;\;8\;\;16\;\;16\;\;26\;\;26\;\;26\;\;200
\\1\;\;1\;\;1\;\;8\;\;26\;\;26\;\;26\;\;48\;\;48\;\;198\;\;200
\\1\;\;1\;\;1\;\;8\;\;8\;\;26\;\;26\;\;26\;\;48\;\;48\;\;198\;\;198\;\;200
\\1\;\;1\;\;1\;\;1\;\;1\;\;8\;\;8\;\;16\;\;16\;\;16\;\;26\;\;26\;\;26\;\;198\;\;200
\\1\;\;1\;\;1\;\;1\;\;1\;\;1\;\;8\;\;8\;\;16\;\;16\;\;16\;\;26\;\;26\;\;26\;\;200
\\1\;\;1\;\;1\;\;1\;\;8\;\;8\;\;16\;\;26\;\;26\;\;26\;\;48\;\;48\;\;198\;\;200
\\1\;\;1\;\;1\;\;8\;\;16\;\;26\;\;26\;\;26\;\;48\;\;48\;\;198\;\;200
\\1\;\;1\;\;1\;\;8\;\;16\;\;26\;\;26\;\;26\;\;200
\\1\;\;1\;\;1\;\;8\;\;8\;\;16\;\;26\;\;198
\\1\;\;1\;\;8\;\;8\;\;16\;\;26\;\;48\;\;198
\\1\;\;1\;\;16\;\;26\;\;48\;\;200
\\1\;\;26\;\;200
\\1\;\;8
\\8
\end{array}
\quad
\begin{array}{c}
16
\\16\;\;26
\\26\;\;48\;\;198
\\26\;\;26\;\;48\;\;198\;\;200
\\1\;\;1\;\;1\;\;16\;\;26\;\;26\;\;200
\\1\;\;1\;\;1\;\;8\;\;8\;\;16\;\;26
\\1\;\;1\;\;8\;\;8\;\;16\;\;16\;\;26\;\;48
\\1\;\;1\;\;16\;\;16\;\;26\;\;26\;\;26\;\;48\;\;200\;\;200
\\1\;\;1\;\;1\;\;26\;\;26\;\;26\;\;48\;\;48\;\;198\;\;200\;\;200
\\1\;\;1\;\;1\;\;8\;\;8\;\;8\;\;26\;\;26\;\;26\;\;48\;\;48\;\;198\;\;198
\\1\;\;1\;\;1\;\;1\;\;1\;\;8\;\;8\;\;8\;\;16\;\;16\;\;16\;\;16\;\;26\;\;26\;\;26\;\;198
\\1\;\;1\;\;1\;\;1\;\;1\;\;1\;\;8\;\;16\;\;16\;\;16\;\;16\;\;26\;\;26\;\;26\;\;26\;\;200
\\1\;\;1\;\;1\;\;8\;\;26\;\;26\;\;26\;\;26\;\;48\;\;48\;\;48\;\;198\;\;200
\\1\;\;1\;\;8\;\;26\;\;26\;\;26\;\;48\;\;48\;\;48\;\;198
\\1\;\;1\;\;8\;\;16\;\;16\;\;16\;\;26\;\;26\;\;26\;\;198
\\1\;\;1\;\;1\;\;8\;\;16\;\;16\;\;16\;\;26\;\;26\;\;198\;\;200
\\1\;\;1\;\;8\;\;16\;\;26\;\;26\;\;48\;\;198\;\;200
\\1\;\;16\;\;26\;\;48\;\;198\;\;200
\\16\;\;26\;\;200
\\16
\end{array}\]

\[\begin{array}{c}
26
\\
1\;\;16\;\;26\;\;48\;\;198
\\1\;\;8\;\;16\;\;26\;\;26\;\;26\;\;48\;\;198
\\1\;\;1\;\;1\;\;8\;\;16\;\;16\;\;26\;\;26\;\;26\;\;48\;\;198
\\1\;\;1\;\;1\;\;8\;\;16\;\;16\;\;26\;\;26\;\;26\;\;26\;\;26\;\;48\;\;198\;\;200
\\1\;\;1\;\;1\;\;1\;\;1\;\;8\;\;16\;\;26\;\;26\;\;26\;\;26\;\;26\;\;48\;\;48\;\;198\;\;200
\\1\;\;1\;\;1\;\;1\;\;1\;\;8\;\;8\;\;8\;\;16\;\;26\;\;26\;\;26\;\;26\;\;48\;\;48\;\;198\;\;200
\\1\;\;1\;\;1\;\;1\;\;1\;\;1\;\;1\;\;8\;\;8\;\;8\;\;16\;\;16\;\;16\;\;26\;\;26\;\;26\;\;26\;\;48\;\;48\;\;200
\\1\;\;1\;\;1\;\;1\;\;1\;\;1\;\;1\;\;8\;\;8\;\;8\;\;16\;\;16\;\;16\;\;26\;\;26\;\;26\;\;26\;\;26\;\;48\;\;48\;\;198\;\;200\;\;200
\\1\;\;1\;\;1\;\;1\;\;1\;\;1\;\;8\;\;8\;\;8\;\;16\;\;16\;\;16\;\;26\;\;26\;\;26\;\;26\;\;26\;\;48\;\;48\;\;48\;\;198\;\;198\;\;200\;\;200
\\1\;\;1\;\;1\;\;1\;\;1\;\;1\;\;1\;\;8\;\;8\;\;8\;\;16\;\;16\;\;16\;\;26\;\;26\;\;26\;\;26\;\;26\;\;26\;\;26\;\;48\;\;48\;\;48\;\;198\;\;198\;\;200\;\;200
\\1\;\;1\;\;1\;\;1\;\;1\;\;1\;\;1\;\;1\;\;8\;\;8\;\;8\;\;16\;\;16\;\;16\;\;16\;\;26\;\;26\;\;26\;\;26\;\;26\;\;26\;\;26\;\;48\;\;48\;\;198\;\;198\;\;200\;\;200
\\1\;\;1\;\;1\;\;1\;\;1\;\;1\;\;1\;\;1\;\;8\;\;8\;\;8\;\;16\;\;16\;\;16\;\;16\;\;26\;\;26\;\;26\;\;26\;\;26\;\;26\;\;48\;\;48\;\;198\;\;198\;\;200
\\1\;\;1\;\;1\;\;1\;\;1\;\;8\;\;8\;\;8\;\;16\;\;16\;\;16\;\;26\;\;26\;\;26\;\;26\;\;26\;\;26\;\;48\;\;48\;\;48\;\;198\;\;198\;\;200
\\1\;\;1\;\;1\;\;1\;\;1\;\;8\;\;16\;\;16\;\;16\;\;26\;\;26\;\;26\;\;26\;\;26\;\;48\;\;48\;\;48\;\;198\;\;200
\\1\;\;1\;\;1\;\;8\;\;16\;\;16\;\;26\;\;26\;\;26\;\;26\;\;26\;\;48\;\;198\;\;198\;\;200
\\1\;\;1\;\;1\;\;8\;\;16\;\;16\;\;26\;\;26\;\;26\;\;48\;\;198\;\;198
\\1\;\;1\;\;8\;\;16\;\;26\;\;26\;\;26\;\;48\;\;198
\\1\;\;16\;\;26\;\;48\;\;198
\\26
\end{array}\]

\[\begin{array}{c}
48
\\26\;\;48
\\1\;\;16\;\;26
\\1\;\;8\;\;16\;\;26
\\1\;\;8\;\;26\;\;48\;\;198
\\1\;\;26\;\;26\;\;48\;\;198\;\;200
\\1\;\;1\;\;1\;\;16\;\;26\;\;26\;\;48\;\;200
\\1\;\;1\;\;1\;\;8\;\;8\;\;16\;\;26\;\;26\;\;48
\\1\;\;1\;\;8\;\;8\;\;16\;\;16\;\;26\;\;26\;\;48
\\1\;\;1\;\;1\;\;16\;\;16\;\;26\;\;26\;\;26\;\;48\;\;200
\\1\;\;1\;\;1\;\;26\;\;26\;\;26\;\;48\;\;48\;\;198\;\;200
\\1\;\;1\;\;8\;\;8\;\;26\;\;26\;\;48\;\;48\;\;198\;\;198
\\1\;\;1\;\;1\;\;8\;\;8\;\;16\;\;16\;\;16\;\;26\;\;26\;\;198
\\1\;\;1\;\;1\;\;1\;\;16\;\;16\;\;16\;\;26\;\;26\;\;26\;\;200
\\1\;\;1\;\;26\;\;26\;\;26\;\;48\;\;198\;\;200
\\1\;\;8\;\;26\;\;48\;\;198
\\1\;\;8\;\;16\;\;26
\\1\;\;16\;\;26\;\;48
\\26\;\;48\;\;48
\\48
\end{array}
\quad\begin{array}{c}
198
\\26\;\;198
\\1\;\;16\;\;26
\\1\;\;16\;\;26
\\26\;\;48
\\1\;\;26\;\;48
\\1\;\;1\;\;26\;\;200
\\1\;\;1\;\;8\;\;198\;\;200
\\1\;\;1\;\;8\;\;8\;\;16\;\;26\;\;198
\\1\;\;1\;\;8\;\;16\;\;16\;\;26\;\;26\;\;200
\\1\;\;1\;\;16\;\;26\;\;26\;\;48\;\;198\;\;200
\\1\;\;8\;\;26\;\;26\;\;48\;\;48\;\;198
\\1\;\;8\;\;16\;\;26\;\;26\;\;48
\\1\;\;1\;\;16\;\;26\;\;26
\\1\;\;8\;\;16\;\;26\;\;48\;\;198
\\1\;\;8\;\;16\;\;26\;\;26\;\;48\;\;198
\\1\;\;16\;\;26\;\;26\;\;198\;\;200
\\1\;\;1\;\;16\;\;26\;\;198\;\;200
\\1\;\;1\;\;26\;\;198
\\198
\end{array}\quad
\begin{array}{c}
200
\\1\;\;200
\\1\;\;8
\\1\;\;8\;\;16
\\1\;\;16\;\;26\;\;200
\\1\;\;26\;\;48\;\;200
\\1\;\;8\;\;26\;\;48\;\;198
\\1\;\;1\;\;1\;\;8\;\;16\;\;16\;\;26\;\;198
\\1\;\;1\;\;1\;\;8\;\;16\;\;16\;\;26\;\;26\;\;200
\\1\;\;1\;\;8\;\;26\;\;26\;\;48\;\;198\;\;200
\\1\;\;1\;\;8\;\;26\;\;26\;\;48\;\;198
\\1\;\;1\;\;8\;\;16\;\;26\;\;26
\\1\;\;1\;\;8\;\;16\;\;26
\\1\;\;8\;\;26\;\;48
\\1\;\;26\;\;48\;\;200
\\1\;\;1\;\;16\;\;26\;\;200
\\1\;\;1\;\;8\;\;16\;\;198\;\;200
\\1\;\;8\;\;16\;\;198\;\;200
\\1\;\;16\;\;200
\\200
\end{array}\]
This is the largest symmetric group where we are able to write the projectives like this, since the projective cover of the trivial module is getting too wide for the page; thus we will just write down the projectives for the alternating groups.

\section{$A_{11}$ and $S_{11}$, characteristic $2$}

We will only give here the projectives for the principal block of the alternating group, leaving as an exercise placing another copy of a number directly above all numbers to generate the projectives for the symmetric group.

There are seven simple modules in the principal block of $S_{11}$ (and $A_{11}$) and they are below.
\begin{center}
\begin{tabular}{ccc}
\hline $\lambda$ &$\dim(D^\lambda)$ & Factors of $S^\lambda$
\\\hline  $(11)$ & $1$ & $1$
\\ $(9,2)$ & $44$ & $44$
\\ $(6,4,1)$ & $144$ & $144,164,186,198,1$
\\ $(7,4)$ & $164$ & $164,1$
\\ $(8,2,1)$ & $186$ & $186,44,1$
\\ $(7,3,1)$ & $198$ & $198,186,164,1^2$
\\ $(5,4,2)$ & $416$ & $416,144,186,198,44,1^2$
\\ \hline
\end{tabular}
\end{center}

The non-principal block has $2$-core $(2,1)$, and has modules as follows for the alternating group.
\begin{center}
\begin{tabular}{ccc}
\hline $\lambda$ &$\dim(D^\lambda)$ & Factors of $S^\lambda$
\\\hline  $(10,1)$ & $10$ & $10$
\\ $(6,5)$ & $16_1\oplus 16_2$ & $32,100$ 
\\ $(8,3)$ & $100$ & $100,10$
\\ $(6,3,2)$ & $848$ & $848,100,32,10$
\\ $(5,3,2,1)$ & $584_1\oplus 584_2$ & $1168,848,32^2,100^2,10^3$
\\ \hline
\end{tabular}
\end{center}

Here are the projectives for the principal $2$-block of the alternating group $A_{11}$, which we can then induce to $S_{11}$, if that is what you need.

\[\begin{array}{c}
1
\\44\;\;164\;\;198\;\;198\;\;416
\\1\;\;1\;\;1\;\;1\;\;1\;\;186\;\;186
\\44\;\;44\;\;144\;\;164\;\;164\;\;198\;\;198\;\;416\;\;416
\\1\;\;1\;\;1\;\;1\;\;1\;\;1\;\;1\;\;1\;\;186\;\;186\;\;186
\\44\;\;44\;\;44\;\;44\;\;144\;\;144\;\;164\;\;164\;\;164\;\;198\;\;198\;\;198\;\;198\;\;416\;\;416
\\1\;\;1\;\;1\;\;1\;\;1\;\;1\;\;1\;\;1\;\;1\;\;1\;\;1\;\;1\;\;186\;\;186\;\;186
\\44\;\;44\;\;44\;\;44\;\;44\;\;144\;\;144\;\;164\;\;164\;\;164\;\;186\;\;198\;\;198\;\;198\;\;198\;\;198\;\;416\;\;416\;\;416
\\1\;\;1\;\;1\;\;1\;\;1\;\;1\;\;1\;\;1\;\;1\;\;1\;\;1\;\;1\;\;186\;\;186\;\;186\;\;186
\\1\;\;44\;\;44\;\;44\;\;44\;\;144\;\;144\;\;164\;\;164\;\;164\;\;186\;\;198\;\;198\;\;198\;\;198\;\;416\;\;416\;\;416
\\1\;\;1\;\;1\;\;1\;\;1\;\;1\;\;1\;\;1\;\;1\;\;1\;\;144\;\;164\;\;186\;\;186\;\;186
\\1\;\;1\;\;44\;\;44\;\;44\;\;44\;\;144\;\;164\;\;164\;\;186\;\;186\;\;198\;\;198\;\;198\;\;198\;\;416\;\;416
\\1\;\;1\;\;1\;\;1\;\;1\;\;1\;\;1\;\;1\;\;144\;\;164\;\;186\;\;186\;\;198
\\1\;\;44\;\;44\;\;44\;\;144\;\;164\;\;164\;\;186\;\;198\;\;198\;\;416\;\;416
\\1\;\;1\;\;1\;\;1\;\;1\;\;1\;\;1\;\;144\;\;164\;\;186
\\44\;\;164\;\;186\;\;198\;\;198\;\;416
\\1
\end{array}\]

\[\begin{array}{c}
44
\\1
\\164\;\;416
\\1\;\;1\;\;186
\\44\;\;44\;\;44\;\;198\;\;198
\\1\;\;1\;\;1\;\;1\;\;186
\\144\;\;144\;\;164\;\;164\;\;416\;\;416
\\1\;\;1\;\;1\;\;1\;\;1\;\;186\;\;186
\\44\;\;44\;\;44\;\;198\;\;198\;\;198\;\;198
\\1\;\;1\;\;1\;\;1\;\;186
\\144\;\;164\;\;164\;\;186\;\;416\;\;416
\\1\;\;1\;\;1\;\;1\;\;186
\\44\;\;44\;\;44\;\;198\;\;198
\\1\;\;1\;\;1\;\;144
\\164\;\;186\;\;416
\\1\;\;44
\\44
\end{array}\quad
\begin{array}{c}
144
\\186
\\198
\\1
\\164\;\;186\;\;416
\\1\;\;1
\\44\;\;44\;\;198
\\1\;\;1\;\;144
\\144\;\;164\;\;186\;\;416
\\1\;\;1\;\;186\;\;198
\\1\;\;44\;\;198
\\1\;\;164\;\;186\;\;186
\\1\;\;144\;\;164\;\;198
\\1\;\;44\;\;186
\\1\;\;144\;\;198
\\186\;\;416
\\144
\end{array}\quad
\begin{array}{c}
164
\\1
\\44\;\;198
\\1\;\;1\;\;186
\\144\;\;164\;\;164\;\;416
\\1\;\;1\;\;1\;\;186
\\44\;\;44\;\;198\;\;198\;\;198
\\1\;\;1\;\;1\;\;186
\\144\;\;164\;\;164\;\;186\;\;416
\\1\;\;1\;\;1\;\;186
\\1\;\;44\;\;44\;\;198\;\;198
\\1\;\;1\;\;144\;\;164\;\;186
\\1\;\;144\;\;164\;\;186\;\;416
\\1\;\;1\;\;186\;\;198
\\1\;\;44\;\;198
\\1\;\;164
\\164
\end{array}\]

\[\begin{array}{c}
186
\\144\;\;198
\\1\;\;1\;\;186
\\44\;\;164\;\;186\;\;198
\\1\;\;1\;\;1\;\;144
\\44\;\;164\;\;186\;\;186\;\;198\;\;416
\\1\;\;1\;\;1\;\;198
\\1\;\;44\;\;44\;\;164\;\;198\;\;416
\\1\;\;1\;\;1\;\;1\;\;144\;\;164\;\;186\;\;186
\\1\;\;44\;\;144\;\;164\;\;186\;\;198\;\;416
\\1\;\;1\;\;1\;\;44\;\;186\;\;186\;\;198\;\;198
\\1\;\;1\;\;44\;\;144\;\;144\;\;164\;\;198
\\1\;\;1\;\;164\;\;186\;\;186\;\;186\;\;416
\\1\;\;144\;\;164\;\;198\;\;198
\\1\;\;44\;\;186\;\;186
\\1\;\;144\;\;198
\\186
\end{array}\quad
\begin{array}{c}
198
\\1\;\;1\;\;186
\\144\;\;164\;\;416
\\1\;\;1\;\;186
\\44\;\;44\;\;198\;\;198\;\;198
\\1\;\;1\;\;1\;\;1\;\;186
\\144\;\;164\;\;164\;\;164\;\;186\;\;416\;\;416
\\1\;\;1\;\;1\;\;1\;\;1\;\;186
\\44\;\;44\;\;44\;\;44\;\;198\;\;198\;\;198\;\;198
\\1\;\;1\;\;1\;\;1\;\;144\;\;186
\\144\;\;164\;\;164\;\;186\;\;186\;\;416\;\;416
\\1\;\;1\;\;1\;\;1\;\;186\;\;198
\\1\;\;44\;\;44\;\;144\;\;198\;\;198
\\1\;\;1\;\;164\;\;186\;\;186
\\144\;\;164\;\;198\;\;416
\\1\;\;1\;\;186
\\198
\end{array}\quad
\begin{array}{c}
416
\\1
\\44\;\;198
\\1\;\;1
\\144\;\;164\;\;416
\\1\;\;1\;\;186
\\44\;\;44\;\;198\;\;198
\\1\;\;1\;\;1\;\;186
\\144\;\;164\;\;416\;\;416
\\1\;\;1\;\;1\;\;186
\\44\;\;44\;\;198\;\;198
\\1\;\;1
\\164\;\;186\;\;416
\\1\;\;1
\\44\;\;198
\\1\;\;144
\\416
\end{array}\]

For the non-principal block of $A_{11}$, we need $\F_4$ for the two pairs of dual modules, of dimension $16$ and $584$. We start by giving the structure of the projectives over $\F_2$ for the modules $10$, $100$, $848$, and the two irreducible but not absolutely irreducible modules of dimension $32$ and $1168$. We do this for ease of translation to the symmetric group case. We then will give the projective covers of $16$ and $584$ over $\F_4$, choosing this labelling so that $\Ext^1(584,16)$ is non-zero (but $\Ext^1(584,16^*)=0$).

\[\begin{array}{c}
10
\\100\;\;100\;\;848
\\10\;\;10\;\;32\;\;32
\\100\;\;100\;\;100\;\;848
\\10\;\;10\;\;32\;\;32
\\100\;\;100\;\;100\;\;848\;\;1168
\\10\;\;10\;\;10\;\;10\;\;10\;\;10\;\;10\;\;32\;\;32\;\;32
\\100\;\;100\;\;100\;\;100\;\;100\;\;100\;\;100\;\;848\;\;848\;\;1168\;\;1168
\\10\;\;10\;\;10\;\;10\;\;10\;\;10\;\;10\;\;10\;\;10\;\;10\;\;10\;\;10\;\;32\;\;32\;\;32\;\;32\;\;32\;\;32
\\100\;\;100\;\;100\;\;100\;\;100\;\;100\;\;100\;\;848\;\;848\;\;1168\;\;1168
\\10\;\;10\;\;10\;\;10\;\;10\;\;10\;\;10\;\;32\;\;32\;\;32
\\100\;\;100\;\;848\;\;1168
\\10
\end{array}\]

\[\begin{array}{c}
32
\\
100\;\;100\;\;848\;\;848
\\
10\;\;10\;\;10\;\;10\;\;32\;\;32\;\;32
\\
100\;\;100\;\;100\;\;100\;\;100\;\;100\;\;848\;\;848
\\
10\;\;10\;\;10\;\;10\;\;32\;\;32\;\;32\;\;32
\\
100\;\;100\;\;100\;\;100\;\;100\;\;100\;\;848\;\;848\;\;1168
\\
10\;\;10\;\;10\;\;10\;\;10\;\;10\;\;32\;\;32\;\;32\;\;32\;\;32
\\
100\;\;100\;\;100\;\;100\;\;100\;\;100\;\;100\;\;100\;\;848\;\;848\;\;1168\;\;1168
\\
10\;\;10\;\;10\;\;10\;\;10\;\;10\;\;10\;\;10\;\;10\;\;10\;\;10\;\;10\;\;32\;\;32\;\;32\;\;32\;\;32\;\;32\;\;32
\\
100\;\;100\;\;100\;\;100\;\;100\;\;100\;\;100\;\;100\;\;848\;\;848\;\;1168\;\;1168
\\
10\;\;10\;\;10\;\;10\;\;10\;\;10\;\;32\;\;32\;\;32\;\;32
\\
100\;\;100\;\;848\;\;848\;\;1168
\\
32
\end{array}\]

\[\begin{array}{c}
100
\\10\;\;10\;\;32
\\100\;\;100\;\;848\;\;848
\\10\;\;10\;\;10\;\;32\;\;32\;\;32
\\100\;\;100\;\;100\;\;100\;\;100\;\;848
\\10\;\;10\;\;10\;\;32\;\;32\;\;32
\\100\;\;100\;\;100\;\;100\;\;100\;\;848\;\;848\;\;1168
\\10\;\;10\;\;10\;\;10\;\;10\;\;10\;\;10\;\;32\;\;32\;\;32\;\;32
\\100\;\;100\;\;100\;\;100\;\;100\;\;100\;\;100\;\;848\;\;1168
\\10\;\;10\;\;10\;\;10\;\;10\;\;10\;\;10\;\;32\;\;32\;\;32\;\;32
\\100\;\;100\;\;100\;\;848\;\;848\;\;1168
\\10\;\;10\;\;32
\\100
\end{array}\quad\begin{array}{c}
848
\\10\;\;32
\\100\;\;100
\\10\;\;32
\\100\;\;848
\\10\;\;32
\\100\;\;100
\\10\;\;10\;\;32
\\100\;\;848\;\;1168
\\10\;\;10\;\;32
\\100\;\;100
\\10\;\;32
\\848
\end{array}\quad\begin{array}{c}
1168
\\10\;\;10\;\;32
\\100\;\;100\;\;1168
\\10\;\;10\;\;10\;\;10\;\;32\;\;32
\\100\;\;100\;\;848\;\;848\;\;1168
\\10\;\;10\;\;10\;\;10\;\;32\;\;32
\\100\;\;100\;\;1168
\\10\;\;10\;\;32
\\1168
\end{array}\]

\[\begin{array}{c}
16
\\100\;\;848
\\10\;\;10\;\;16\;\;16^*\;\;16^*
\\100\;\;100\;\;100\;\;848
\\10\;\;10\;\;16\;\;16\;\;16^*\;\;16^*
\\100\;\;100\;\;100\;\;584^*\;\;848
\\10\;\;10\;\;10\;\;16\;\;16\;\;16^*\;\;16^*\;\;16^*
\\100\;\;100\;\;100\;\;100\;\;584\;\;584\;\;848
\\10\;\;10\;\;10\;\;10\;\;10\;\;10\;\;16\;\;16\;\;16\;\;16\;\;16\;\;16^*\;\;16^*
\\100\;\;100\;\;100\;\;100\;\;584^*\;\;584^*\;\;848
\\10\;\;10\;\;10\;\;16\;\;16^*\;\;16^*\;\;16^*
\\100\;\;584\;\;848
\\16
\end{array}\quad\begin{array}{c}
584
\\10\;\;16
\\100\;\;584^*
\\10\;\;10\;\;16^*\;\;16^*
\\100\;\;584\;\;848
\\10\;\;10\;\;16\;\;16
\\100\;\;584^*
\\10\;\;16^*
\\584
\end{array}\]
\section{$A_{12}$ and $S_{12}$, characteristic $2$}

This group is getting towards the edge of what can be easily constructed using today's computers. What makes it even more annoying is that there are simple modules in the principal $2$-block of $A_{12}$ that are not realizable over $\F_2$.

There are eleven simple modules in the principal block of $S_{12}$, three of which split into two dual modules upon restriction to $A_{12}$, and they are below.
\begin{center}
\begin{tabular}{cccc}
\hline $\lambda$ &$\dim(D^\lambda)$ & Factors of $S^\lambda$ & $\dim(P(D^\lambda))$
\\\hline  $(12)$ & $1$ & $1$ & $204288$
\\ $(11,1)$ & $10$ & $10$ & $159232$
\\ $(7,5)$ & $32=16\oplus 16^*$ & $32,164,100,1$ & $145408$
\\ $(10,2)$ & $44$ & $44,10$ & $69120$
\\ $(9,3)$ & $100$ & $100,44,10$ & $116224$
\\ $(8,4)$ & $164$ & $164,100,10,1$ & $59904$
\\ $(6,5,1)$ & $288=144\oplus144^*$ & $288,570,164,100,32,1$& $55296$
\\ $(6,4,2)$ & $416$ & $416,1046,570,288,164,100,44,32,10,1^3$ & $50688$
\\ $(8,3,1)$ & $570$ & $570,164,100,44,10,1^3$ & $48640$
\\ $(7,3,2)$ & $1046$ & $1046,570,164,100,32,10,1^3$ & $45056$
\\ $(5,4,2,1)$ & $2368=1184\oplus 1184^*$ & $2368,1046,570,416^2,288,164,100^3,44^2,32^2,10^5,1^5$ & $45056$
\\ \hline
\end{tabular}
\end{center}

Each of the dual pairs amalgamate over $\F_2$. We will not give the structures of the projectives for the $1$-, $10$-, $16$- and $100$-dimensional simple modules because it takes too long to find their socle layers. However, we do give the $\Ext^1$ matrix below.

The non-principal block not of defect $0$ has $2$-core $(3,2,1)$, and has three modules as follows for the alternating and symmetric groups.
\begin{center}
\begin{tabular}{cccc}
\hline $\lambda$ &$\dim(D^\lambda)$ & Factors of $S^\lambda$ & $\dim(P(D^\lambda))$
\\\hline $(9,2,1)$ & $320$ & $320$ & $7680$
\\ $(7,4,1)$ & $1408$ & $1408$ & $6656$
\\ $(5,4,3)$ & $1792$ & $1792,320$ & $5632$
\\ \hline
\end{tabular}
\end{center}

There is also a block with defect $1$ for $S_{12}$ and of defect zero for $A_{12}$, with simple of dimension $5632$ and given by the partition $(6,3,2,1)$.

Here is the $\Ext^1$-matrix for the principal block for $p=2$.

\begin{center}
\begin{tabular}{c|cccccccccccccc}
&$1$&$10$&$16$&$16^*$&$44$&$100$&
$144$&$144^*$&$164$&$416$&$570$&$1046$&$1184$&$1184^*$
\\\hline $1$ &$0$&$1$&$0$&$0$&$1$&$0$&$0$&$0$&$1$&$1$&$2$&$1$&$0$&$0$
\\ $10$ &$1$&$0$&$0$&$0$&$1$&$2$&$0$&$0$&$0$&$0$&$0$&$1$&$1$&$1$
\\ $16$&$0$&$0$&$0$&$1$&$0$&$1$&$1$&$0$&$0$&$0$&$0$&$1$&$1$&$0$
\\ $16^*$&$0$&$0$&$1$&$0$&$0$&$1$&$0$&$1$&$0$&$0$&$0$&$1$&$0$&$1$
\\ $44$&$1$&$1$&$0$&$0$&$0$&$0$&$0$&$0$&$0$&$0$&$0$&$0$&$0$&$0$
\\ $100$&$0$&$2$&$1$&$1$&$0$&$0$&$0$&$0$&$1$&$0$&$0$&$0$&$0$&$0$
\\ $144$&$0$&$0$&$1$&$0$&$0$&$0$&$0$&$0$&$0$&$1$&$1$&$0$&$0$&$0$
\\ $144^*$&$0$&$0$&$0$&$1$&$0$&$0$&$0$&$0$&$0$&$1$&$1$&$0$&$0$&$0$
\\ $164$&$1$&$0$&$0$&$0$&$0$&$1$&$0$&$0$&$0$&$0$&$0$&$0$&$0$&$0$
\\ $416$&$1$&$0$&$0$&$0$&$0$&$0$&$1$&$1$&$0$&$0$&$0$&$0$&$1$&$1$
\\ $570$&$2$&$0$&$0$&$0$&$0$&$0$&$1$&$1$&$0$&$0$&$0$&$0$&$0$&$0$
\\ $1046$&$1$&$1$&$1$&$1$&$0$&$0$&$0$&$0$&$0$&$0$&$0$&$0$&$0$&$0$
\\ $1184$ &$0$&$1$&$1$&$0$&$0$&$0$&$0$&$0$&$0$&$1$&$0$&$0$&$0$&$0$
\\ $1184^*$&$0$&$1$&$0$&$1$&$0$&$0$&$0$&$0$&$0$&$1$&$0$&$0$&$0$&$0$
\end{tabular}
\end{center}

%\begin{footnotesize}
%\[\begin{array}{c}
%10\;\;44\;\;164\;\;416\;\;570\;\;570\;\;1046
%\\1
%\end{array}\]
%\end{footnotesize}

We now give the socle layers of the projectives that we have constructed.

\begin{footnotesize}
\[\begin{array}{c}
44
\\1
\\10\;\;164\;\;416\;\;570
\\1\;\;1\;\;1\;\;1\;\;100\;\;100
\\10\;\;32\;\;44\;\;44\;\;164\;\;570
\\1\;\;1\;\;1\;\;100\;\;288\;\;1046\;\;1046
\\10\;\;10\;\;32\;\;32\;\;44\;\;44\;\;164\;\;164\;\;416\;\;570\;\;570
\\1\;\;1\;\;1\;\;1\;\;1\;\;1\;\;1\;\;1\;\;10\;\;10\;\;100\;\;100\;\;100\;\;288
\\10\;\;10\;\;32\;\;44\;\;44\;\;44\;\;100\;\;164\;\;164\;\;164\;\;416\;\;416\;\;570\;\;570\;\;570\;\;2368
\\1\;\;1\;\;1\;\;1\;\;1\;\;1\;\;1\;\;1\;\;1\;\;1\;\;10\;\;10\;\;10\;\;10\;\;32\;\;100\;\;100\;\;100\;\;100\;\;288\;\;1046
\\10\;\;10\;\;32\;\;32\;\;44\;\;44\;\;44\;\;44\;\;44\;\;100\;\;100\;\;100\;\;164\;\;164\;\;164\;\;416\;\;416\;\;570\;\;570\;\;1046
\\1\;\;1\;\;1\;\;1\;\;1\;\;1\;\;1\;\;1\;\;1\;\;1\;\;1\;\;1\;\;10\;\;10\;\;10\;\;10\;\;10\;\;10\;\;32\;\;100\;\;100\;\;100\;\;288\;\;288\;\;1046\;\;1046
\\10\;\;10\;\;10\;\;10\;\;32\;\;32\;\;32\;\;32\;\;44\;\;44\;\;44\;\;44\;\;44\;\;44\;\;44\;\;100\;\;164\;\;164\;\;164\;\;164\;\;416\;\;416\;\;416\;\;570\;\;570\;\;570\;\
\;570\;\;1046\;\;1046\;\;2368\;\;2368
\\1\;\;1\;\;1\;\;1\;\;1\;\;1\;\;1\;\;1\;\;1\;\;1\;\;1\;\;1\;\;1\;\;1\;\;1\;\;1\;\;1\;\;1\;\;10\;\;10\;\;10\;\;10\;\;10\;\;10\;\;10\;\;10\;\;32\;\;32\;\;100\;\;100\;\;100\;\;100\;\;100\;\;100\;\;288
\\10\;\;10\;\;32\;\;32\;\;44\;\;44\;\;44\;\;44\;\;100\;\;100\;\;100\;\;100\;\;100\;\;100\;\;164\;\;164\;\;164\;\;164\;\;164\;\;416\;\;416\;\;416\;\;416\;\;416\;\;570\;\
\;570\;\;570\;\;1046
\\1\;\;1\;\;1\;\;1\;\;1\;\;1\;\;1\;\;1\;\;1\;\;1\;\;1\;\;1\;\;1\;\;1\;\;10\;\;10\;\;10\;\;10\;\;10\;\;10\;\;10\;\;10\;\;32\;\;32\;\;100\;\;100\;\;100\;\;100\;\;288\;\;288\;\;288\;\;1046\;\;1046
\\10\;\;10\;\;10\;\;32\;\;32\;\;32\;\;44\;\;44\;\;44\;\;44\;\;44\;\;44\;\;44\;\;100\;\;164\;\;164\;\;164\;\;416\;\;570\;\;570\;\;570\;\;1046\;\;1046\;\;1046\;\;2368\;\;2368
\\1\;\;1\;\;1\;\;1\;\;1\;\;1\;\;1\;\;1\;\;1\;\;1\;\;1\;\;1\;\;1\;\;10\;\;10\;\;10\;\;10\;\;10\;\;10\;\;10\;\;32\;\;100\;\;100\;\;100\;\;100\;\;288
\\10\;\;10\;\;32\;\;32\;\;44\;\;44\;\;44\;\;100\;\;100\;\;100\;\;100\;\;164\;\;164\;\;164\;\;164\;\;416\;\;416\;\;416\;\;416\;\;570\;\;570\;\;570
\\1\;\;1\;\;1\;\;1\;\;1\;\;1\;\;1\;\;1\;\;10\;\;10\;\;10\;\;10\;\;10\;\;32\;\;100\;\;100\;\;288\;\;1046
\\10\;\;32\;\;44\;\;44\;\;44\;\;44\;\;100\;\;164\;\;570\;\;1046\;\;2368
\\1\;\;1\;\;1\;\;1\;\;10\;\;10\;\;10\;\;100\;\;288
\\44\;\;100\;\;164\;\;416\;\;570
\\1\;\;10
\\44
\end{array}\]
\end{footnotesize}

%\[
%\begin{array}{c}
%\\ 10,\;\;10,\;\;32,\;\;164
%\\ 100
%\end{array}\]
%
%[ 1, 1, 1, 2, 2, 2, 2, 2, 2, 2, 2, 3, 3, 3, 3, 3, 4, 6, 6, 9 ]
%[ 1, 4, 5, 5, 5, 5, 10, 10, 11 ]
%[ 2, 2, 3, 6 ]
%[ 5 ]

\[
\begin{array}{c}
144
\\416\;\;570
\\1\;\;1
\\16\;\;44\;\;164
\\1\;\;1\;\;100\;\;144\;\;144^*\;\;1046
\\10\;\;16\;\;16^*\;\;44\;\;164\;\;416\;\;570
\\1\;\;1\;\;1\;\;1\;\;100\;\;100\;\;144^*
\\10\;\;16\;\;16^*\;\;44\;\;164\;\;416\;\;570\;\;570
\\1\;\;1\;\;1\;\;1\;\;100\;\;144\;\;144^*\;\;1046
\\10\;\;10\;\;16\;\;16^*\;\;16^*\;\;44\;\;44\;\;164\;\;164\;\;570\;\;1184^*
\\1\;\;1\;\;1\;\;1\;\;1\;\;10\;\;100\;\;100\;\;144\;\;144^*\;\;1046
\\10\;\;16\;\;16^*\;\;44\;\;100\;\;100\;\;164\;\;164\;\;416\;\;416\;\;570\;\;570
\\1\;\;1\;\;1\;\;1\;\;1\;\;1\;\;10\;\;10\;\;10\;\;16\;\;16^*\;\;100\;\;100
\\10\;\;16\;\;16\;\;16^*\;\;44\;\;44\;\;44\;\;164\;\;570\;\;1046\;\;1184
\\1\;\;1\;\;1\;\;1\;\;10\;\;10\;\;100\;\;100\;\;144\;\;144^*\;\;1046
\\10\;\;16\;\;16^*\;\;44\;\;100\;\;164\;\;164\;\;416\;\;570
\\1\;\;1\;\;1\;\;1\;\;1\;\;10\;\;100\;\;100\;\;144\;\;144
\\10\;\;16\;\;16^*\;\;44\;\;164\;\;416\;\;570\;\;570
\\1\;\;1\;\;10\;\;16^*\;\;100\;\;144\;\;144^*
\\16\;\;44\;\;164\;\;416\;\;1046
\\1\;\;1\;\;144\;\;144^*\;\;144^*\;\;1184^*
\\16\;\;416\;\;570
\\144
\end{array}\]

\[\begin{array}{c}
164
\\1\;\;100
\\10\;\;32\;\;44\;\;164\;\;570
\\1\;\;1\;\;1\;\;100\;\;1046
\\10\;\;32\;\;44\;\;100\;\;164\;\;416\;\;570
\\1\;\;1\;\;1\;\;1\;\;1\;\;1\;\;10\;\;32\;\;100\;\;288
\\10\;\;44\;\;44\;\;100\;\;164\;\;164\;\;416\;\;570\;\;570\;\;1046
\\1\;\;1\;\;1\;\;1\;\;1\;\;1\;\;10\;\;32\;\;100\;\;100\;\;288\;\;1046
\\10\;\;10\;\;10\;\;32\;\;32\;\;32\;\;44\;\;44\;\;44\;\;100\;\;164\;\;164\;\;164\;\;416\;\;570\;\;570\;\;1046
\\1\;\;1\;\;1\;\;1\;\;1\;\;1\;\;1\;\;1\;\;1\;\;10\;\;10\;\;100\;\;100\;\;100\;\;100\;\;100\;\;288\;\;1046
\\10\;\;10\;\;10\;\;32\;\;32\;\;32\;\;44\;\;44\;\;44\;\;100\;\;164\;\;164\;\;164\;\;164\;\;416\;\;416\;\;570\;\;570\;\;570\;\;2368
\\1\;\;1\;\;1\;\;1\;\;1\;\;1\;\;1\;\;1\;\;1\;\;1\;\;1\;\;1\;\;1\;\;10\;\;10\;\;10\;\;10\;\;32\;\;100\;\;100\;\;100\;\;100\;\;288\;\;288\;\;1046
\\10\;\;10\;\;32\;\;32\;\;44\;\;44\;\;44\;\;44\;\;100\;\;100\;\;100\;\;164\;\;164\;\;164\;\;164\;\;164\;\;416\;\;416\;\;570\;\;570\;\;570\;\;570\;\;1046
\\1\;\;1\;\;1\;\;1\;\;1\;\;1\;\;1\;\;1\;\;1\;\;1\;\;1\;\;1\;\;1\;\;10\;\;10\;\;10\;\;10\;\;32\;\;100\;\;100\;\;100\;\;100\;\;100\;\;288\;\;288\;\;1046\;\;10\
46
\\10\;\;10\;\;10\;\;10\;\;32\;\;32\;\;32\;\;32\;\;44\;\;44\;\;44\;\;44\;\;44\;\;100\;\;164\;\;164\;\;164\;\;164\;\;164\;\;416\;\;416\;\;570\;\;570\;\;570\;\;1046\;\;2368
\\1\;\;1\;\;1\;\;1\;\;1\;\;1\;\;1\;\;1\;\;1\;\;1\;\;1\;\;1\;\;1\;\;10\;\;10\;\;10\;\;10\;\;10\;\;10\;\;32\;\;32\;\;100\;\;100\;\;100\;\;100\;\;100\;\;288\;\;1046
\\10\;\;10\;\;32\;\;32\;\;44\;\;44\;\;44\;\;100\;\;100\;\;100\;\;100\;\;164\;\;164\;\;164\;\;164\;\;416\;\;416\;\;416\;\;570\;\;570\;\;570\;\;1046
\\1\;\;1\;\;1\;\;1\;\;1\;\;1\;\;1\;\;1\;\;1\;\;1\;\;1\;\;10\;\;10\;\;10\;\;10\;\;32\;\;100\;\;100\;\;100\;\;100\;\;288\;\;288\;\;1046
\\10\;\;10\;\;10\;\;32\;\;32\;\;32\;\;44\;\;44\;\;44\;\;44\;\;164\;\;164\;\;164\;\;164\;\;416\;\;570\;\;570\;\;1046\;\;2368
\\1\;\;1\;\;1\;\;1\;\;1\;\;1\;\;1\;\;1\;\;10\;\;10\;\;100\;\;100\;\;100\;\;100\;\;288\;\;1046
\\10\;\;10\;\;32\;\;32\;\;44\;\;100\;\;164\;\;164\;\;416\;\;570\;\;570
\\1\;\;1\;\;1\;\;1\;\;10\;\;100\;\;100\;\;288\;\;1046
\\10\;\;32\;\;44\;\;164\;\;570
\\1\;\;100
\\164
\end{array}\]

\[\begin{array}{c}
416
\\1
\\44
\\1\;\;288\;\;1046
\\10\;\;32\;\;164\;\;416\;\;570
\\1\;\;1\;\;1\;\;1\;\;100\;\;100
\\10\;\;32\;\;44\;\;164\;\;416\;\;570\;\;2368
\\1\;\;1\;\;1\;\;10\;\;10\;\;100\;\;288\;\;1046
\\10\;\;32\;\;44\;\;44\;\;100\;\;100\;\;164\;\;416\;\;570
\\1\;\;1\;\;1\;\;1\;\;1\;\;10\;\;10\;\;10\;\;32\;\;100\;\;288\;\;1046
\\10\;\;10\;\;32\;\;32\;\;44\;\;44\;\;100\;\;164\;\;164\;\;416\;\;416\;\;570\;\;570\;\;1046\;\;2368
\\1\;\;1\;\;1\;\;1\;\;1\;\;1\;\;1\;\;1\;\;1\;\;10\;\;10\;\;10\;\;10\;\;32\;\;100\;\;100\;\;100
\\10\;\;32\;\;44\;\;44\;\;44\;\;100\;\;100\;\;100\;\;164\;\;164\;\;416\;\;416\;\;416\;\;570\;\;1046
\\1\;\;1\;\;1\;\;1\;\;1\;\;1\;\;1\;\;1\;\;10\;\;10\;\;10\;\;10\;\;10\;\;32\;\;100\;\;100\;\;288\;\;288\;\;1046
\\10\;\;10\;\;32\;\;32\;\;44\;\;44\;\;44\;\;44\;\;44\;\;164\;\;164\;\;416\;\;570\;\;570\;\;1046\;\;1046\;\;2368\;\;2368
\\1\;\;1\;\;1\;\;1\;\;1\;\;1\;\;1\;\;1\;\;1\;\;1\;\;10\;\;10\;\;10\;\;10\;\;10\;\;32\;\;100\;\;100\;\;100
\\10\;\;32\;\;44\;\;100\;\;100\;\;100\;\;100\;\;164\;\;164\;\;164\;\;416\;\;416\;\;416\;\;570\;\;570
\\1\;\;1\;\;1\;\;1\;\;1\;\;1\;\;1\;\;10\;\;10\;\;10\;\;10\;\;32\;\;100\;\;100\;\;288\;\;1046
\\10\;\;32\;\;44\;\;44\;\;44\;\;44\;\;164\;\;570\;\;1046\;\;1046\;\;2368
\\1\;\;1\;\;1\;\;1\;\;10\;\;10\;\;10\;\;10\;\;32\;\;100\;\;288
\\100\;\;100\;\;100\;\;164\;\;416\;\;416\;\;570
\\1\;\;1\;\;1\;\;10\;\;10\;\;32\;\;288
\\10\;\;32\;\;44\;\;416\;\;416\;\;1046
\\1\;\;288\;\;2368
\\416
\end{array}\]

\[\begin{array}{c}
570
\\1\;\;1
\\44\;\;164
\\1\;\;1\;\;100\;\;288\;\;1046
\\10\;\;32\;\;44\;\;164\;\;416\;\;570\;\;570
\\1\;\;1\;\;1\;\;1\;\;1\;\;1\;\;100
\\10\;\;32\;\;44\;\;44\;\;100\;\;164\;\;164\;\;416\;\;570\;\;570
\\1\;\;1\;\;1\;\;1\;\;1\;\;1\;\;10\;\;32\;\;100\;\;288\;\;1046
\\10\;\;10\;\;32\;\;44\;\;44\;\;44\;\;164\;\;164\;\;416\;\;570\;\;1046
\\1\;\;1\;\;1\;\;1\;\;1\;\;1\;\;1\;\;100\;\;100\;\;100\;\;100\;\;288\;\;288\;\;1046
\\10\;\;10\;\;10\;\;32\;\;32\;\;32\;\;44\;\;44\;\;100\;\;164\;\;164\;\;164\;\;416\;\;416\;\;570\;\;570\;\;570\;\;570
\\1\;\;1\;\;1\;\;1\;\;1\;\;1\;\;1\;\;1\;\;1\;\;1\;\;1\;\;10\;\;10\;\;100\;\;100\;\;100\;\;100\;\;288\;\;1046
\\10\;\;10\;\;10\;\;32\;\;32\;\;32\;\;44\;\;44\;\;44\;\;44\;\;100\;\;164\;\;164\;\;164\;\;164\;\;416\;\;570\;\;570\;\;570\;\;2368
\\1\;\;1\;\;1\;\;1\;\;1\;\;1\;\;1\;\;1\;\;1\;\;10\;\;10\;\;10\;\;10\;\;32\;\;100\;\;100\;\;100\;\;288\;\;288\;\;1046\;\;1046
\\10\;\;10\;\;32\;\;32\;\;44\;\;44\;\;44\;\;100\;\;100\;\;100\;\;164\;\;164\;\;164\;\;416\;\;416\;\;570\;\;570\;\;570\;\;1046
\\1\;\;1\;\;1\;\;1\;\;1\;\;1\;\;1\;\;1\;\;1\;\;1\;\;10\;\;10\;\;10\;\;10\;\;32\;\;100\;\;100\;\;100\;\;100\;\;288
\\10\;\;10\;\;32\;\;32\;\;44\;\;44\;\;44\;\;100\;\;164\;\;164\;\;164\;\;416\;\;416\;\;570\;\;570\;\;1046\;\;2368
\\1\;\;1\;\;1\;\;1\;\;1\;\;1\;\;1\;\;10\;\;10\;\;10\;\;32\;\;100\;\;100\;\;288\;\;1046
\\10\;\;32\;\;44\;\;44\;\;44\;\;100\;\;164\;\;164\;\;416\;\;570\;\;570\;\;1046
\\1\;\;1\;\;1\;\;1\;\;1\;\;1\;\;1\;\;10\;\;100\;\;100\;\;100\;\;288\;\;288
\\10\;\;32\;\;44\;\;164\;\;164\;\;416\;\;570\;\;570
\\1\;\;1\;\;1\;\;10\;\;32\;\;100
\\44\;\;164\;\;570\;\;1046
\\1\;\;1\;\;288
\\570
\end{array}\]

\[\begin{array}{c}
1046
\\10\;\;32
\\1\;\;1\;\;100
\\100\;\;164\;\;416\;\;570
\\1\;\;1\;\;10\;\;32
\\44\;\;44\;\;100\;\;1046
\\1\;\;1\;\;10\;\;32\;\;164\;\;288\;\;1046
\\1\;\;1\;\;10\;\;10\;\;32\;\;100\;\;164\;\;416\;\;570
\\1\;\;1\;\;1\;\;1\;\;100\;\;100\;\;100\;\;164\;\;570\;\;2368
\\1\;\;1\;\;10\;\;10\;\;10\;\;10\;\;32\;\;32\;\;44\;\;164\;\;416\;\;570
\\1\;\;1\;\;1\;\;44\;\;100\;\;100\;\;100\;\;288\;\;416\;\;1046\;\;1046
\\1\;\;1\;\;10\;\;10\;\;10\;\;10\;\;32\;\;32\;\;44\;\;44\;\;100\;\;164\;\;570
\\1\;\;1\;\;1\;\;1\;\;10\;\;32\;\;44\;\;44\;\;100\;\;164\;\;288\;\;416\;\;1046\;\;1046\;\;2368
\\1\;\;1\;\;1\;\;1\;\;1\;\;10\;\;10\;\;10\;\;10\;\;32\;\;100\;\;100\;\;164\;\;164\;\;416\;\;570\;\;570
\\1\;\;1\;\;1\;\;1\;\;10\;\;32\;\;44\;\;100\;\;100\;\;100\;\;100\;\;100\;\;164\;\;416\;\;416\;\;570
\\1\;\;1\;\;1\;\;1\;\;10\;\;10\;\;10\;\;10\;\;10\;\;32\;\;32\;\;44\;\;44\;\;100\;\;164\;\;288\;\;1046
\\1\;\;1\;\;10\;\;32\;\;44\;\;44\;\;44\;\;100\;\;164\;\;288\;\;570\;\;1046\;\;1046\;\;1046\;\;2368
\\1\;\;1\;\;1\;\;1\;\;10\;\;10\;\;10\;\;10\;\;32\;\;32\;\;100\;\;164\;\;416\;\;570
\\1\;\;1\;\;1\;\;100\;\;100\;\;100\;\;100\;\;100\;\;164\;\;416\;\;416\;\;570
\\1\;\;1\;\;1\;\;10\;\;10\;\;10\;\;10\;\;32\;\;32\;\;44\;\;164
\\1\;\;44\;\;100\;\;1046\;\;1046\;\;2368
\\1\;\;10\;\;10\;\;32\;\;164\;\;288
\\1\;\;100\;\;100\;\;416\;\;570
\\1\;\;10\;\;32
\\1046
\end{array}\]

%\[\begin{array}{c}
%2368
%\\10\;\;10
%\\100\;\;100
%\\10\;\;10
%\\2368
%\\10\;\;10
%\\100\;\;100\;\;416\;\;416
%\\1\;\;1\;\;10\;\;10\;\;10\;\;10\;\;32\;\;32
%\\44\;\;44\;\;100\;\;100\;\;1046\;\;1046\;\;2368
%\\1\;\;1\;\;1\;\;1\;\;10\;\;10\;\;10\;\;10\;\;32
%\\100\;\;100\;\;164\;\;164\;\;416\;\;416\;\;2368
%\\1\;\;1\;\;1\;\;1\;\;10\;\;10\;\;10\;\;10\;\;10\;\;10\;\;32\;\;32\;\;288
%\\44\;\;44\;\;44\;\;44\;\;100\;\;100\;\;100\;\;100\;\;570\;\;570\;\;1046\;\;1046\;\;2368
%\\1\;\;1\;\;1\;\;1\;\;1\;\;1\;\;10\;\;10\;\;10\;\;10\;\;10\;\;10\;\;32
%\\100\;\;100\;\;164\;\;164\;\;416\;\;416\;\;416\;\;416\;\;2368
%\\1\;\;1\;\;1\;\;1\;\;1\;\;1\;\;10\;\;10\;\;10\;\;10\;\;10\;\;10\;\;32\;\;32\;\;32\;\;288
%\\44\;\;44\;\;44\;\;44\;\;100\;\;100\;\;100\;\;100\;\;570\;\;570\;\;1046\;\;1046\;\;2368
%\\1\;\;1\;\;1\;\;1\;\;10\;\;10\;\;10\;\;10\;\;10\;\;10\;\;32
%\\100\;\;100\;\;164\;\;164\;\;416\;\;416\;\;2368
%\\1\;\;1\;\;1\;\;1\;\;10\;\;10\;\;10\;\;10\;\;32\;\;32
%\\44\;\;44\;\;100\;\;100\;\;1046\;\;1046\;\;2368
%\\10\;\;10\;\;10\;\;10\;\;32
%\\1\;\;1\;\;100\;\;100\;\;288
%\\10\;\;10\;\;32\;\;416\;\;416
%\\2368
%\end{array}\]

\[\quad\begin{array}{c}
1184
\\10
\\100
\\10
\\1184^*
\\10
\\100\;\;416
\\1\;\;10\;\;10\;\;16\;\;16
\\44\;\;100\;\;1046\;\;1184
\\1\;\;1\;\;10\;\;10\;\;16^*
\\100\;\;164\;\;416\;\;1184
\\1\;\;1\;\;10\;\;10\;\;10\;\;16^*\;\;16^*\;\;144
\\44\;\;44\;\;100\;\;100\;\;570\;\;1046\;\;1184^*
\\1\;\;1\;\;1\;\;10\;\;10\;\;10\;\;16
\\100\;\;164\;\;416\;\;416\;\;1184^*
\\1\;\;1\;\;1\;\;10\;\;10\;\;10\;\;16\;\;16\;\;16\;\;144^*
\\44\;\;44\;\;100\;\;100\;\;570\;\;1046\;\;1184
\\1\;\;1\;\;10\;\;10\;\;10\;\;16^*
\\100\;\;164\;\;416\;\;1184
\\1\;\;1\;\;10\;\;10\;\;16^*\;\;16^*
\\44\;\;100\;\;1046\;\;1184^*
\\10\;\;10\;\;16
\\1\;\;100\;\;144
\\10\;\;16\;\;416
\\1184
\end{array}\]

We also have the non-principal block, with projectives as follows:

\[\begin{array}{c}
320
\\1408\;\;1792
\\320\;\;320
\\1408\;\;1792
\\320\end{array}
\quad\begin{array}{c}
1408
\\320
\\1792
\\320\;\;1408
\\1408\end{array}
\quad\begin{array}{c}
1792
\\320
\\1408
\\320
\\1792
\end{array}\]

The time it took to construct the projectives on my computer is as follows. Note that there are preliminary things like constructing the Cartan matrix and finding good subgroups for condensation that, if one constructs several projectives in one instance of Magma, will not be repeated, so use these numbers as a guide only.

\begin{center}
\begin{tabular}{cc}
\hline $D^\lambda$ & time to construct $P(D^\lambda)$ (seconds)
\\\hline $1$ & 64295
\\ $10$ & 9196
\\ $16$ (over $\F_4$) & 27314
\\ $32$ & 30555
\\ $44$ & 9750
\\ $100$ & 14141
\\ $164$ & 6570
\\ $144$ (over $\F_4$) & 2676
\\ $288$ & 785
\\ $416$ & 7097
\\ $570$ & 2328
\\ $1046$& 7647
\\ $1184$ (over $\F_4$) & 2490
\\ $2368$ & 618
\\ \hline
\end{tabular}
\end{center}

\section{$A_{10}$, characteristic $3$}

In $A_{10}$ in characteristic $3$ there are two blocks: the principal block with $3$-core $(1)$, of defect $4$, and the non-principal block with $3$-core $(3,1)$.

\begin{center}
\begin{tabular}{ccc|ccc}
\hline $\lambda$ &$\dim(D^\lambda)$& Factors of $S^\lambda$&$\lambda$ & $\dim(D^\lambda)$ & Factors of $S^\lambda$
\\\hline $(10)$&$1$&$1$ & $(9,1)$ & $9$ & $9$
\\ $(8,2)$ & $34$& $34,1$ & $(6,4)$ & $90$ & $90$
\\ $(7,3)$ & $41$& $41,34$ &$(6,2^2)$& $126$ & $126,90,9$
\\ $(7,2,1)$ & $84$& $84,41,34,1$&$(4^2,2)$ & $36$ & $36,126,90$
\\ $(6,2,1^2)$ & $224$& $224,84,41,1$ & $(3^2,2,1^2)$ & $279$ &$279,126,36,9$
\\\hline
\end{tabular}
\end{center}

The non-principal block is Morita equivalent to the two non-principal blocks of $S_{10}$, which are labelled by the $3$-cores $(3,1)$ and its conjugate $(2,1^2)$. These are Morita equivalent to the principal $3$-block of $S_8$, by standard Scopes moves \cite{scopes1991}.

The $\Ext^1$-matrices are easy to describe, and we do this now.

\begin{center}\begin{tabular}{c|ccccc}
Block 1& $1$&$34$&$41$&$84$&$224$
\\\hline $1$&$0$&$1$&$1$&$1$&$1$
\\ $34$ & $1$&$0$&$1$&$0$&$1$
\\ $41$&$1$&$1$&$0$&$1$&$1$
\\ $84$&$1$&$0$&$1$&$1$&$0$
\\ $224$&$1$&$1$&$1$&$0$&$1$
\end{tabular}\end{center}

\begin{center}\begin{tabular}{c|ccccc}
Block 2& $9$&$36$&$90$&$126$&$279$
\\\hline $9$&$0$&$0$&$0$&$1$&$1$
\\ $36$&$0$&$0$&$0$&$1$&$1$
\\ $90$&$0$&$0$&$0$&$1$&$0$
\\ $126$&$1$&$1$&$1$&$0$&$0$
\\ $279$&$1$&$1$&$0$&$0$&$0$
\end{tabular}\end{center}

Finally, here are the radical and socle layers of the projectives, which coincide in this case.

\[ 
\begin{array}{c}
1 
\\34\;\;41\;\;84\;\;224
\\1\;\;1\;\;1\;\;1\;\;34\;\;41\;\;41\;\;84\;\;224\;\;224
\\1\;\;1\;\;1\;\;34\;\;34\;\;34\;\;41\;\;41\;\;41\;\;84\;\;84\;\;224\;\;224
\\1\;\;1\;\;1\;\;1\;\;34\;\;41\;\;41\;\;84\;\;224\;\;224
\\34\;\;41\;\;84\;\;224
\\1\end{array}\qquad \begin{array}{c}
34
\\1\;\;41\;\;224
\\1\;\;34\;\;34\;\;41\;\;84\;\;224
\\1\;\;1\;\;1\;\;34\;\;41\;\;41\;\;84\;\;224
\\1\;\;34\;\;34\;\;41\;\;84\;\;224
\\1\;\;41\;\;224
\\34
\end{array}\]

\[ 
\begin{array}{c}
41
\\1\;\;34\;\;84\;\;224
\\1\;\;1\;\;34\;\;41\;\;41\;\;41\;\;84\;\;224
\\1\;\;1\;\;1\;\;34\;\;34\;\;41\;\;84\;\;84\;\;224\;\;224
\\1\;\;1\;\;34\;\;41\;\;41\;\;41\;\;84\;\;224
\\1\;\;34\;\;84\;\;224
\\41
\end{array}\qquad
\begin{array}{c}
84
\\1\;\;41\;\;84
\\1\;\;34\;\;41\;\;84\;\;84\;\;224
\\1\;\;1\;\;34\;\;41\;\;41\;\;84\;\;224
\\1\;\;34\;\;41\;\;84\;\;84\;\;224
\\1\;\;41\;\;84
\\84
\end{array}\quad \begin{array}{c}
224
\\1\;\;34\;\;41\;\;224
\\1\;\;1\;\;34\;\;41\;\;84\;\;224\;\;224
\\1\;\;1\;\;34\;\;41\;\;41\;\;84\;\;224\;\;224
\\1\;\;1\;\;34\;\;41\;\;84\;\;224\;\;224
\\1\;\;34\;\;41\;\;224
\\224
\end{array}\]

\[\begin{array}{c}9
\\126\;\;279
\\9\;\;9\;\;36\;\;90
\\126\;\;279
\\9\end{array}\quad
\begin{array}{c}36
\\126\;\;279
\\9\;\;36\;\;36\;\;90
\\126\;\;279
\\36\end{array}\quad
\begin{array}{c}90
\\126
\\9\;\;36\;\;90
\\126
\\90\end{array}\quad
\begin{array}{c}126
\\9\;\;36\;\;90
\\126\;\;126\;\;279
\\9\;\;36\;\;90
\\126\end{array}\quad
\begin{array}{c}279
\\9\;\;36
\\126\;\;279
\\9\;\;36
\\279
\end{array}\]

These were described for the symmetric group by Tan in \cite{tan2000}, and restriction is simply a case of replacing the module for $S_{10}$ by its restriction.

\section{$A_{11}$, characteristic $3$}

In $A_{11}$ in characteristic $3$ there are two blocks, the principal block with $3$-core $(2)$, of defect $4$, the non-principal block with $3$-core $(3,1,1)$ of defect $2$, and a block with defect $1$ and $3$-core $(4,2,1^2)$.

The principal block of $A_{11}$ is Morita equivalent to that of $S_{11}$, and we label the modules from this block by $3$-regular partitions with $3$-core $(2)$, rather than $(1^2)$, so that the $10$-dimensional simple module has label $(5^2,1)$, rather than $(10,1)$.

\begin{center}
\begin{tabular}{ccc}
\hline $\lambda$ &$\dim(D^\lambda)$& Factors of $S^\lambda$
\\\hline $(11)$&$1$&$1$ 
\\ $(8,3)$ &$109$ & $109,1$
\\ $(8,2,1)$ &$120$ & $120,109$
\\ $(7,3,1)$ & $320$ & $320,120,109,1$
\\$(6,3,1^2)$ & $791$ & $791,320,120,1$
\\ $(5^2,1)$ & $10$ & $10,320$
\\ $(5,4,1^2)$ & $34$ & $34,10,791,320,1^2$
\\ $(5,3^2)$ & $210$ & $210,10,320,120$
\\ $(5,3,2,1)$ & $714$ & $714,210,34,10,791,320,120,109,1^2$
\\ $(4,3,2^2)$ & $131$ & $131,714,210,34,120,109,1^2$
\\\hline
\end{tabular}
\end{center}
%$(9,1^2)$ gives $45$
%$(6,4,1)$ gives $693$
%$(6,3,2)$ gives $126\oplus 126^*$
%
%$(7,2,1^2)$ gives $594$

We now give the $\Ext^1$-matrix, for those readers who do not want to look at the projectives directly.
\begin{center}\begin{tabular}{c|cccccccccc}
Block 1& $1$&$10$&$34$&$109$&$120$&$131$&$210$&$320$&$714$&$791$
\\\hline $1$&$0$&$0$&$0$&$1$&$1$&$1$&$0$&$0$&$1$&$1$
\\$10$&$0$&$0$&$1$&$0$&$0$&$0$&$1$&$1$&$0$&$0$
\\$34$&$0$&$1$&$0$&$0$&$0$&$1$&$0$&$0$&$1$&$1$
\\$109$&$1$&$0$&$0$&$0$&$0$&$0$&$0$&$0$&$0$&$0$
\\$120$&$1$&$0$&$0$&$0$&$0$&$0$&$1$&$1$&$0$&$0$
\\$131$&$1$&$0$&$1$&$0$&$0$&$0$&$1$&$0$&$0$&$0$
\\$210$&$0$&$1$&$0$&$0$&$1$&$1$&$0$&$0$&$1$&$0$
\\$320$&$0$&$1$&$0$&$0$&$1$&$0$&$0$&$0$&$0$&$1$
\\$714$&$1$&$0$&$1$&$0$&$0$&$0$&$1$&$0$&$0$&$0$
\\$791$&$1$&$0$&$1$&$0$&$0$&$0$&$0$&$1$&$0$&$0$
\end{tabular}\end{center}

\begin{center}\begin{tabular}{c|cccc}
Block 2& $45$&$126$&$126^*$&$693$
\\\hline $45$&$0$&$1$&$1$&$0$
\\ $126$&$1$&$0$&$0$&$1$
\\ $126^*$&$1$&$0$&$0$&$1$
\\ $693$&$0$&$1$&$1$&$0$
\end{tabular}\end{center}

For those who do, here they are.
\[ 
\begin{array}{c}
\\1
\\109\;\;120\;\;714
\\1\;\;1\;\;1
\\109\;\;109\;\;120\;\;714
\\1\;\;1\;\;1\;\;1\;\;34\;\;210\;\;320
\\109\;\;109\;\;120\;\;120\;\;131\;\;714\;\;791
\\1\;\;1\;\;1\;\;1\;\;1\;\;34\;\;34\;\;210\;\;210\;\;320
\\10\;\;109\;\;109\;\;120\;\;120\;\;131\;\;714\;\;714\;\;791
\\1\;\;1\;\;1\;\;1\;\;1\;\;34\;\;210\;\;320
\\10\;\;10\;\;109\;\;109\;\;120\;\;120\;\;131\;\;131\;\;714\;\;714\;\;791
\\1\;\;1\;\;1\;\;1\;\;1\;\;34\;\;34\;\;210\;\;210\;\;320
\\109\;\;120\;\;131\;\;714\;\;791
\\1\end{array}
\quad
\begin{array}{c}10
\\34\;\;210\;\;320
\\120\;\;714
\\1
\\10\;\;10\;\;10\;\;109\;\;131\;\;791
\\1\;\;1\;\;34\;\;34\;\;34\;\;210\;\;210\;\;320\;\;320
\\10\;\;10\;\;10\;\;120\;\;131\;\;714\;\;791
\\34\;\;210\;\;320
\\10
\end{array}\]

\[\begin{array}{c}
34
\\714
\\1
\\10\;\;109\;\;131\;\;791
\\1\;\;1\;\;34\;\;210\;\;320
\\120\;\;714
\\1\;\;34\;\;34\;\;210
\\10\;\;10\;\;10\;\;109\;\;120\;\;131\;\;131\;\;714\;\;791
\\1\;\;1\;\;34\;\;34\;\;34\;\;210\;\;210\;\;320
\\10\;\;131\;\;714\;\;791
\\34
\end{array}\quad
\begin{array}{c}
109
\\1
\\120\;\;714
\\1\;\;1
\\109\;\;109
\\1\;\;1\;\;34\;\;210\;\;320
\\120\;\;120\;\;714
\\1\;\;1
\\10\;\;109\;\;109\;\;131\;\;791
\\1\;\;1\;\;34\;\;210\;\;320
\\120\;\;714
\\1
\\109
\end{array}
\quad
\begin{array}{c}
120
\\1
\\109
\\1\;\;210\;\;320
\\120\;\;120\;\;714
\\1\;\;1
\\10\;\;109\;\;109\;\;131\;\;791
\\1\;\;1\;\;34\;\;210\;\;320
\\120\;\;120\;\;714
\\1\;\;1\;\;34\;\;210\;\;320
\\10\;\;109\;\;120\;\;120\;\;131\;\;791
\\1\;\;210\;\;320
\\120
\end{array}\]
\[\begin{array}{c}
131
\\1\;\;34\;\;210
\\120\;\;714
\\1
\\10\;\;109\;\;131\;\;131
\\1\;\;1\;\;34\;\;34\;\;210\;\;210
\\10\;\;120\;\;131\;\;131\;\;714
\\1\;\;34\;\;210
\\131\end{array}
\quad
\begin{array}{c}
210
\\120\;\;714
\\1
\\10\;\;109\;\;131
\\1\;\;1\;\;34\;\;210\;\;320
\\120\;\;714
\\1\;\;34\;\;210\;\;210
\\10\;\;10\;\;109\;\;120\;\;131\;\;131\;\;714\;\;791
\\1\;\;1\;\;34\;\;34\;\;210\;\;210\;\;210\;\;320
\\10\;\;120\;\;131\;\;714
\\210\end{array}
\quad
\begin{array}{c}
320
\\120
\\1
\\10\;\;109\;\;791
\\1\;\;34\;\;210\;\;320
\\120\;\;714
\\1\;\;320
\\10\;\;10\;\;109\;\;120\;\;791
\\1\;\;34\;\;210\;\;320\;\;320
\\10\;\;120\;\;791
\\320\end{array}\]

\[\begin{array}{c}
714
\\1
\\109
\\1\;\;34\;\;210
\\120\;\;714
\\1
\\10\;\;109\;\;131\;\;791
\\1\;\;1\;\;34\;\;210\;\;320
\\120\;\;714\;\;714
\\1\;\;1\;\;34\;\;210
\\10\;\;109\;\;131\;\;714\;\;791
\\1\;\;34\;\;210
\\714
\end{array}\quad
\begin{array}{c}
791
\\1\;\;34\;\;320
\\120\;\;714
\\1
\\10\;\;109\;\;791
\\1\;\;34\;\;210\;\;320
\\10\;\;120\;\;714\;\;791
\\1\;\;34\;\;320
\\791
\end{array}\]

And the non-principal block with defect $2$.
\[
\begin{array}{c}
45
\\126\;\;126^*
\\45\;\;693\;\;693
\\126\;\;126^*
\\45
\end{array}\quad
\begin{array}{c}
126
\\45\;\;693
\\126\;\;126^*\;\;126^*
\\45\;\;693
\\126
\end{array}\quad
\begin{array}{c}
126^*
\\45\;\;693
\\126\;\;126\;\;126^*
\\45\;\;693
\\126^*
\end{array}\quad
\begin{array}{c}
693
\\126\;\;126^*
\\45\;\;45\;\;693
\\126\;\;126^*
\\693
\end{array}\]

We should note that three of these, of dimensions $109$, $120$ and $131$, were determined by Tan in \cite{tan2000a}. 

\section{$A_{12}$, characteristic $3$}

And we are back up to the largest case you can feel comfortable with using current computers, and here we do not produce the socle layers of the projective cover of the trivial module.

There are four blocks for $A_{12}$ in characteristic $3$, of defects $5$, $2$, $1$ and $0$. The block of defect zero has a module of dimension $2673$, and the block of defect $1$ has two modules of dimensions $891$ and $3564$.

\begin{center}
\begin{tabular}{ccc}
\hline $\lambda$ &$\dim(D^\lambda)$& Factors of $S^\lambda$
\\\hline $(12)$&$1$&$1$ 
\\ $(11,1)$ &$10$&$10,1$
\\ $(10,1^2)$ & $45$ & $45,10$
\\ $(9,3)$&$143$&$143,10,1$
\\ $(9,2,1)$&$120$&$120,143,45,10,1^2$
%\\ $(9,1^3)$ & $120$ & $120,45$
\\ $(8,4)$ & $131$ & $131,143,1$
\\ $(8,2^2)$ &$210$&$210,131,120,143,10,1^2$
\\ $(7,4,1)$ & $1013$ & $1013,131,120,143,1$
\\ $(7,3,2)$ & $126\oplus 126^*$ & $252,1013,210,131,120,143,45,10,1$
%\\ $(6^2)$ & $1$ & $1,131$
%\\ $(6,5,1)$ & $1155-$
\\ $(6,4,1^2)$ & $1936$ & $1936,1013,120,10,1$
%\\ $(6,3^2)$ & $1650-$
\\ $(6,3,2,1)$ & $714\oplus 714^*$ & $1428,1936,252,1013,210^2,131,120^2,143,45,10^2,1^4$
\\ \hline
\end{tabular}
\end{center}

The principal block has thirteen simple modules. Here we give the $\Ext^1$ matrix for these modules, since the projectives take several pages to write out, and we do not include $P(1)$, although of course this only means that $\Ext^1(1,1)$ is not given here, and this is ero for all groups $G$ with $O^p(G)=G$.

\begin{center}\begin{tabular}{c|ccccccccccccc}
Block 1& $1$&$10$&$45$&$120$&$126$&$126^*$&$131$&$143$&$210$&$714$&$714^*$&$1013$&$1936$
\\\hline $1$&$0$&$1$&$0$&$1$&$0$&$0$&$1$&$1$&$0$&$1$&$1$&$0$&$1$
\\$10$&$1$&$0$&$1$&$0$&$0$&$0$&$0$&$0$&$1$&$0$&$0$&$1$&$1$
\\$45$&$0$&$1$&$0$&$1$&$1$&$1$&$0$&$0$&$0$&$0$&$0$&$0$&$0$
\\$120$&$1$&$0$&$1$&$0$&$0$&$0$&$0$&$0$&$2$&$0$&$0$&$1$&$0$
\\$126$&$0$&$0$&$1$&$0$&$0$&$0$&$0$&$0$&$1$&$0$&$0$&$1$&$0$
\\$126^*$&$0$&$0$&$1$&$0$&$0$&$0$&$0$&$0$&$1$&$0$&$0$&$1$&$0$
\\$131$&$1$&$0$&$0$&$0$&$0$&$0$&$0$&$1$&$1$&$0$&$0$&$1$&$0$
\\$143$&$1$&$0$&$0$&$0$&$0$&$0$&$1$&$0$&$0$&$0$&$0$&$0$&$0$
\\$210$&$0$&$1$&$0$&$2$&$1$&$1$&$1$&$0$&$0$&$1$&$1$&$0$&$0$
\\$714$&$1$&$0$&$0$&$0$&$0$&$0$&$0$&$0$&$1$&$0$&$0$&$0$&$0$
\\$714^*$&$1$&$0$&$0$&$0$&$0$&$0$&$0$&$0$&$1$&$0$&$0$&$0$&$0$
\\$1013$&$0$&$1$&$0$&$1$&$1$&$1$&$1$&$0$&$0$&$0$&$0$&$0$&$0$
\\$1936$&$1$&$1$&$0$&$0$&$0$&$0$&$0$&$0$&$0$&$0$&$0$&$0$&$0$
\end{tabular}\end{center}

We now give the socle layers of the projective modules in the principal block except for the projective cover of the trivial module.
\[\begin{array}{c}
10
\\1\;\;1936
\\1\;\;45\;\;143\;\;210
\\1\;\;10\;\;10\;\;120\;\;120\;\;714\;\;714^*\;\;1013
\\1\;\;1\;\;1\;\;10\;\;120\;\;120\;\;131\;\;714\;\;714^*
\\1\;\;1\;\;1\;\;10\;\;126\;\;126^*\;\;131\;\;143\;\;143\;\;143\;\;1936
\\1\;\;1\;\;1\;\;1\;\;45\;\;126\;\;126^*\;\;143\;\;143\;\;210\;\;210\;\;210\;\;1013\;\;1013\;\;1936\;\;1936
\\1\;\;1\;\;1\;\;10\;\;10\;\;45\;\;45\;\;120\;\;120\;\;120\;\;143\;\;210\;\;210\;\;210\;\;714\;\;714\;\;714^*\;\;714^*\;\;1013
\\1\;\;1\;\;1\;\;1\;\;1\;\;10\;\;10\;\;10\;\;10\;\;120\;\;120\;\;120\;\;131\;\;131\;\;714\;\;714\;\;714^*\;\;714^*
\\1\;\;1\;\;1\;\;10\;\;10\;\;120\;\;126\;\;126^*\;\;131\;\;143\;\;143\;\;143\;\;1936\;\;1936
\\1\;\;1\;\;1\;\;1\;\;1\;\;1\;\;45\;\;45\;\;126\;\;126^*\;\;143\;\;143\;\;210\;\;210\;\;210\;\;1013\;\;1936
\\1\;\;1\;\;10\;\;10\;\;10\;\;10\;\;10\;\;120\;\;120\;\;120\;\;126\;\;126\;\;126^*\;\;126^*\;\;131\;\;131\;\;143\;\;210\;\;210\;\;714\;\;714^*\;\;1013
\\1\;\;1\;\;1\;\;1\;\;45\;\;45\;\;45\;\;120\;\;210\;\;210\;\;210\;\;714\;\;714^*\;\;1013\;\;1013\;\;1013
\\1\;\;10\;\;10\;\;10\;\;10\;\;120\;\;120\;\;126\;\;126\;\;126^*\;\;126^*\;\;131\;\;131\;\;143
\\1\;\;45\;\;210\;\;1013\;\;1936
\\10\end{array}\]

\[\begin{array}{c}
45
\\10\;\;120
\\1
\\126\;\;126^*\;\;143
\\1\;\;45\;\;210\;\;210\;\;1013\;\;1013
\\10\;\;120\;\;120\;\;120\;\;714\;\;714^*
\\1\;\;1\;\;1\;\;10\;\;10\;\;131\;\;131
\\126\;\;126^*\;\;143\;\;143\;\;143\;\;1936\;\;1936
\\1\;\;1\;\;1\;\;45\;\;45\;\;210\;\;210
\\10\;\;10\;\;120\;\;120\;\;126\;\;126^*\;\;714\;\;714^*
\\1\;\;1\;\;1\;\;45\;\;45\;\;210\;\;210\;\;1013\;\;1013
\\10\;\;10\;\;10\;\;120\;\;120\;\;120\;\;126\;\;126\;\;126^*\;\;126^*\;\;131\;\;131\;\;143
\\1\;\;45\;\;45\;\;210\;\;210\;\;1013\;\;1013
\\10\;\;120\;\;126\;\;126^*
\\45\end{array}\]

\[\begin{array}{c}
120
\\1 
\\143\;\;1936 
\\1\;\;1\;\;45\;\;210\;\;210\;\;1013 
\\10\;\;10\;\;120\;\;120\;\;120\;\;120\;\;120\;\;714\;\;714\;\;714^*\;\;714^*
\\1\;\;1\;\;1\;\;1\;\;1\;\;1\;\;10\;\;10\;\;131\;\;131\;\;131 
\\126\;\;126\;\;126^*\;\;126^*\;\;143\;\;143\;\;143\;\;143\;\;143\;\;1936\;\;1936\;\;1936 
\\1\;\;1\;\;1\;\;1\;\;1\;\;1\;\;1\;\;1\;\;45\;\;45\;\;45\;\;210\;\;210\;\;210\;\;210\;\;210\;\;210\;\;210\;\;1013,
1013\;\;1013 
\\10\;\;10\;\;10\;\;120\;\;120\;\;120\;\;120\;\;120\;\;120\;\;120\;\;143\;\;714\;\;714\;\;714\;\;714\;\;714^*\;\;714^*\;\;714^*\;\;714^*
\\1\;\;1\;\;1\;\;1\;\;1\;\;1\;\;1\;\;1\;\;10\;\;10\;\;10\;\;131\;\;131\;\;131\;\;131 
\\10\;\;120\;\;120\;\;126\;\;126\;\;126^*\;\;126^*\;\;143\;\;143\;\;143\;\;143\;\;143\;\;1936\;\;1936\;\;1936 
\\1\;\;1\;\;1\;\;1\;\;1\;\;1\;\;1\;\;45\;\;45\;\;210\;\;210\;\;210\;\;210\;\;210\;\;210\;\;210\;\;1013\;\;1013 
\\10\;\;10\;\;10\;\;120\;\;120\;\;120\;\;120\;\;120\;\;120\;\;126\;\;126\;\;126^*\;\;126^*\;\;131\;\;131\;\;143\;\;714\;\;714\;\;714^*\;\;714^*
\\1\;\;1\;\;1\;\;1\;\;1\;\;10\;\;45\;\;45\;\;45\;\;131\;\;210\;\;210\;\;210\;\;210\;\;1013\;\;1013\;\;1013 
\\10\;\;10\;\;120\;\;120\;\;120\;\;120\;\;126\;\;126\;\;126^*\;\;126^*\;\;131\;\;131\;\;143\;\;1936 
\\1\;\;45\;\;210\;\;210\;\;1013 
\\120\end{array}\]

\[\begin{array}{c}
126
\\45\;\;210\;\;1013
\\10\;\;120\;\;120\;\;714
\\1\;\;1\;\;10\;\;131
\\126\;\;126^*\;\;126^*\;\;143\;\;143\;\;1936
\\1\;\;1\;\;1\;\;45\;\;210\;\;210\;\;210\;\;1013
\\10\;\;120\;\;120\;\;126^*\;\;714\;\;714^*\;\;714^*
\\1\;\;1\;\;10\;\;45\;\;131\;\;210\;\;1013
\\10\;\;10\;\;120\;\;120\;\;126\;\;126\;\;126\;\;126\;\;126^*\;\;131\;\;143\;\;143\;\;1936
\\1\;\;1\;\;1\;\;45\;\;45\;\;210\;\;210\;\;210\;\;210\;\;1013\;\;1013
\\10\;\;10\;\;120\;\;120\;\;126\;\;126^*\;\;126^*\;\;126^*\;\;131\;\;714
\\45\;\;210\;\;1013
\\126\end{array}\]

\[\begin{array}{c}
131
\\143
\\1\;\;1\;\;210\;\;1013
\\10\;\;120\;\;120\;\;120\;\;714\;\;714^*
\\1\;\;1\;\;1\;\;1\;\;10\;\;131\;\;131
\\126\;\;126^*\;\;143\;\;143\;\;143\;\;143\;\;1936
\\1\;\;1\;\;1\;\;1\;\;1\;\;45\;\;45\;\;210\;\;210\;\;210\;\;210\;\;1013
\\10\;\;10\;\;120\;\;120\;\;120\;\;120\;\;714\;\;714\;\;714^*\;\;714^*
\\1\;\;1\;\;1\;\;1\;\;10\;\;131\;\;131\;\;131\;\;131
\\1\;\;126\;\;126^*\;\;143\;\;143\;\;143\;\;210\;\;1013\;\;1936
\\1\;\;1\;\;1\;\;10\;\;10\;\;120\;\;120\;\;126\;\;126^*\;\;131\;\;131\;\;143\;\;210\;\;210\;\;1013
\\1\;\;1\;\;1\;\;45\;\;45\;\;120\;\;210\;\;210\;\;210\;\;714\;\;714^*\;\;1013
\\1\;\;10\;\;10\;\;120\;\;120\;\;126\;\;126^*\;\;131\;\;131\;\;131
\\1\;\;143\;\;210\;\;1013
\\131\end{array}\]

\[\begin{array}{c}
143
\\1
\\120\;\;714\;\;714^*
\\1\;\;1\;\;1\;\;10\;\;131
\\143\;\;143\;\;143\;\;1936\;\;1936
\\1\;\;1\;\;1\;\;1\;\;1\;\;45\;\;210\;\;210\;\;210\;\;1013
\\10\;\;10\;\;10\;\;120\;\;120\;\;120\;\;120\;\;120\;\;143\;\;714\;\;714\;\;714\;\;714\;\;714\;\;714
\\1\;\;1\;\;1\;\;1\;\;1\;\;1\;\;1\;\;1\;\;10\;\;10\;\;131\;\;131\;\;131\;\;131
\\10\;\;120\;\;126\;\;126\;\;126^*\;\;126^*\;\;143\;\;143\;\;143\;\;143\;\;143\;\;143\;\;1936\;\;1936
\\1\;\;1\;\;1\;\;1\;\;1\;\;1\;\;1\;\;1\;\;45\;\;45\;\;45\;\;210\;\;210\;\;210\;\;210\;\;210\;\;210\;\;210\;\;1013,
1013
\\10\;\;10\;\;10\;\;120\;\;120\;\;120\;\;120\;\;120\;\;143\;\;714\;\;714\;\;714\;\;714^*\;\;714^*\;\;714^*
\\1\;\;1\;\;1\;\;1\;\;1\;\;1\;\;10\;\;10\;\;131\;\;131\;\;131\;\;1013
\\10\;\;120\;\;126\;\;126\;\;126^*\;\;126^*\;\;131\;\;143\;\;143\;\;143\;\;1936\;\;1936
\\1\;\;1\;\;1\;\;1\;\;45\;\;210\;\;210\;\;210\;\;1013
\\10\;\;120\;\;143\;\;714\;\;714^*
\\1\;\;131
\\143\end{array}\]

\[\begin{array}{c}
210
\\10\;\;120\;\;120\;\;714\;\;714^*
\\1\;\;1\;\;1\;\;131
\\126\;\;126^*\;\;143\;\;143\;\;143\;\;1936\;\;1936
\\1\;\;1\;\;1\;\;1\;\;1\;\;1\;\;45\;\;45\;\;210\;\;210\;\;210\;\;210\;\;210\;\;1013\;\;1013
\\10\;\;10\;\;10\;\;120\;\;120\;\;120\;\;120\;\;120\;\;120\;\;120\;\;714\;\;714\;\;714\;\;714\;\;714^*\;\;714^*\;\;714^*\;\;714^*
\\1\;\;1\;\;1\;\;1\;\;1\;\;1\;\;1\;\;1\;\;10\;\;10\;\;10\;\;131\;\;131\;\;131\;\;131
\\126\;\;126\;\;126\;\;126^*\;\;126^*\;\;126^*\;\;143\;\;143\;\;143\;\;143\;\;143\;\;143\;\;143\;\;1936\;\;1936,
1936
\\1\;\;1\;\;1\;\;1\;\;1\;\;1\;\;1\;\;1\;\;1\;\;45\;\;45\;\;210\;\;210\;\;210\;\;210\;\;210\;\;210\;\;210\;\;210\;\;210\;\;1013\;\;1013
\\10\;\;10\;\;10\;\;120\;\;120\;\;120\;\;120\;\;120\;\;120\;\;120\;\;126\;\;126\;\;131\;\;714\;\;714\;\;714\;\;714\;\;714^*\;\;714^*\;\;714^*\;\;714^*
\\1\;\;1\;\;1\;\;1\;\;1\;\;1\;\;1\;\;10\;\;10\;\;45\;\;45\;\;131\;\;131\;\;210\;\;210\;\;210\;\;210\;\;1013\;\;1013
\\10\;\;10\;\;10\;\;120\;\;120\;\;120\;\;120\;\;126\;\;126\;\;126\;\;126\;\;126^*\;\;126^*\;\;126^*\;\;126^*\;\;131\;\;131\;\;131\;\;143\;\;143\;\;143\;\;1936\;\;1936
\\1\;\;1\;\;1\;\;1\;\;45\;\;45\;\;210\;\;210\;\;210\;\;210\;\;210\;\;210\;\;1013\;\;1013
\\10\;\;120\;\;120\;\;126\;\;126^*\;\;131\;\;714\;\; 714^*
\\210\end{array}\]

\[\begin{array}{c}
714
\\1
\\143\;\;1936
\\1\;\;1\;\;210
\\10\;\;120\;\;120\;\;714\;\;714\;\;714^*\;\;714^*
\\1\;\;1\;\;1\;\;1\;\;10\;\;131
\\126\;\;143\;\;143\;\;143\;\;1936\;\;1936
\\1\;\;1\;\;1\;\;1\;\;1\;\;45\;\;210\;\;210\;\;210\;\;210\;\;1013
\\10\;\;10\;\;120\;\;120\;\;120\;\;120\;\;714\;\;714\;\;714\;\;714^*\;\;714^*\;\;714^*
\\1\;\;1\;\;1\;\;1\;\;1\;\;10\;\;10\;\;131\;\;131
\\126\;\;126^*\;\;126^*\;\;143\;\;143\;\;143\;\;1936\;\;1936
\\1\;\;1\;\;1\;\;1\;\;45\;\;210\;\;210\;\;210\;\;210\;\;1013
\\10\;\;120\;\;120\;\;714\;\;714\;\;714^*\;\;714^*
\\1\;\;1\;\;10\;\;131
\\126\;\;143\;\;1936
\\1\;\;210
\\714\end{array}\]

\[\begin{array}{c}
1013
\\120
\\1\;\;10\;\;131
\\126\;\;126^*\;\;143\;\;1936
\\1\;\;1\;\;45\;\;45\;\;210\;\;210\;\;1013
\\10\;\;10\;\;120\;\;120\;\;120\;\;714\;\;714^*
\\1\;\;1\;\;1\;\;10\;\;131
\\126\;\;126^*\;\;143\;\;143\;\;1936
\\1\;\;1\;\;210\;\;210\;\;1013\;\;1013
\\10\;\;120\;\;120\;\;126\;\;126^*\;\;131\;\;143\;\;714\;\;714^*
\\1\;\;1\;\;1\;\;10\;\;45\;\;45\;\;131\;\;210\;\;210\;\;1013\;\;1013
\\10\;\;10\;\;10\;\;120\;\;120\;\;120\;\;126\;\;126\;\;126^*\;\;126^*\;\;131\;\;143\;\;1936
\\1\;\;1\;\;45\;\;45\;\;210\;\;210\;\;1013\;\;1013
\\10\;\;120\;\;126\;\;126^*\;\;131
\\1013\end{array}\]

\[\begin{array}{c}
1936
\\1
\\10\;\;120\;\;714\;\;714^*
\\1\;\;1
\\143\;\;143\;\;1936
\\1\;\;1\;\;1\;\;210\;\;210\;\;1013
\\10\;\;120\;\;120\;\;120\;\;714\;\;714\;\;714^*\;\;714^*
\\1\;\;1\;\;1\;\;1\;\;10\;\;10\;\;131
\\126\;\;126^*\;\;143\;\;143\;\;1936\;\;1936\;\;1936
\\1\;\;1\;\;1\;\;1\;\;45\;\;45\;\;210\;\;210\;\;210\;\;1013
\\10\;\;10\;\;120\;\;120\;\;120\;\;714\;\;714\;\;714^*\;\;714^*
\\1\;\;1\;\;1\;\;10\;\;131
\\126\;\;126^*\;\;143\;\;143\;\;1936
\\1\;\;1\;\;210\;\;210\;\;1013
\\120\;\;714\;\;714^*
\\1\;\;10
\\1936\end{array}\]

The non-principal block of defect $2$ is much smaller of course, and the projectives are as follows.
\[\begin{array}{c}
54
\\1431\;\;1728
\\54\;\;54\;\;297\;\;945
\\1431\;\;1728
\\54
\end{array}\quad
\begin{array}{c}
297
\\1728
\\54\;\;297\;\;945
\\1728
\\297
\end{array}\quad
\begin{array}{c}
945
\\1431\;\;1728
\\54\;\;297\;\;945\;\;945
\\1431\;\;1728
\\945
\end{array}\quad
\begin{array}{c}
1431
\\54\;\;945
\\1431\;\;1728
\\54\;\;945
\\1431
\end{array}\quad
\begin{array}{c}
1728
\\54\;\;945\;\;297
\\1431\;\;1728\;\;1728
\\54\;\;297\;\;945
\\1728
\end{array}
\]

\section{$A_{10}$, characteristic $5$}

Since these can all be understood from the work of Scopes \cite{scopes1995} we simply give the structure of the projectives and the dimensions of the $D^\lambda$, with the factors of $S^\lambda$ for the symmetric group. Where the module for the symmetric group splits into two upon restriction to $A_n$, we note this.
\begin{center}
\begin{tabular}{ccc}
\hline $\lambda$ &$\dim(D^\lambda)$& Factors of $S^\lambda$
\\\hline $(10)$&$1$&$1$ 
\\ $(9,1)$ &$8$&$8,1$
\\ $(8,1^2)$ & $28$ & $28,8$
\\ $(7,1^3)$ & $56$ & $56,28$
\\ $(6,1^4)$ &  $35_1\oplus 35_2$ &  $70,56$
\\ $(5^2)$ & $34$ & $34,8$
\\ $(5,4,1)$ & $217$     & $34,28,8,1$
\\ $(5,3,1^2)$ & $133_1\oplus 133_2$ & $266,217,56,28$
%\\ $(5,2,1^3)$
%\\ $(4^2,2)$
%\\ $(4,3,2,1)$
%\\ $(4,2^2,1^2)$
%\\ $(3^2,2^3)$
%\\ $(3,2^3,1)$
%\\\hline
%
%\\ $(8,2)$ & $35$ & $35$
%\\ $(7,3)$
%\\ $(7,2,1)$
\\ \hline
\end{tabular}
\end{center}

\[
\begin{array}{c}1
\\8\;\;217
\\1\;\;1\;\;28\;\;34\;\;34
\\8\;\;217
\\1
\end{array}\quad
\begin{array}{c}
8
\\1\;\;28\;\;34
\\8\;\;8\;\;217
\\1\;\;28\;\;34
\\8
\end{array}\quad
\begin{array}{c}
28
\\8\;\;56\;\;217
\\1\;\;28\;\;28\;\;34\;\;133_1\;\;133_2
\\8\;\;56\;\;217
\\28
\end{array}\quad
\begin{array}{c}
34
\\8\;\;217\;\;217
\\1\;\;1\;\;28\;\;34\;\;34\;\;34\;\;133_1\;\;133_2
\\8\;\;217\;\;217
\\34
\end{array}\quad
\begin{array}{c}
35_i
\\56
\\35_{3-i}\;\;133_i
\\56
\\35_i
\end{array}\]
\[\begin{array}{c}
56
\\28\;\;35_1\;\;35_2\;\;133_1\;\;133_2
\\56\;\;56\;\;56\;\;217
\\28\;\;35_1\;\;35_2\;\;133_1\;\;133_2
\\56
\end{array}\quad
\begin{array}{c}
133_i
\\56\;\;217
\\28\;\;34\;\;35_i\;\;133_1\;\;133_2
\\56\;\;217
\\133_i
\end{array}\quad
\begin{array}{c}
217
\\1\;\;28\;\;34\;\;34\;\;133_1\;\;133_2
\\8\;\;56\;\;217\;\;217\;\;217
\\1\;\;28\;\;34\;\;34\;\;133_1\;\;133_2
\\217
\end{array}\]

\section{$A_{11}$, characteristic $5$}

Here two modules for $S_{11}$ from the principal block do not remain irreducible when restricted. One restricts to two self-dual modules and one restricts to two dual modules, as we see below.

\begin{center}
\begin{tabular}{ccc}
\hline $\lambda$ &$\dim(D^\lambda)$& Factors of $S^\lambda$
\\\hline $(11)$&$1$&$1$ 
\\ $(9,2)$ & $43$ & $43,1$
\\ $(8,2,1)$ & $188$ & $188,43$
\\ $(7,2,1^2)$ & $406$ & $406,188$
\\ $(6,5)$ & $89$ & $89,43$
\\ $(6,4,1)$ & $372$ & $372,89,188,43,1$
\\ $(6,3,1^2)$ & $133_1\oplus 133_2$ & $266,372,406,188$
\\ $(6,2,1^3)$ & $126\oplus 126^*$ & $252,266,406$
\\ \hline
\end{tabular}
\end{center}

\[
\begin{array}{c}
1
\\43\;\;372
\\1\;\;1\;\;89\;\;89\;\;188
\\43\;\;372
\\1
\end{array}\quad
\begin{array}{c}
43
\\1\;\;89\;\;188
\\43\;\;43\;\;372
\\1\;\;89\;\;188
\\43
\end{array}\quad
\begin{array}{c}
89
\\43\;\;372\;\;372
\\1\;\;1\;\;89\;\;89\;\;89\;\;133_2\;\;133_1\;\;188
\\43\;\;372\;\;372
\\89
\end{array}\quad
\begin{array}{c}
126
\\133_2
\\126^*\;\;406
\\133_1
\\126
\end{array}\quad
\begin{array}{c}
133_1
\\126\;\;372\;\;406
\\89\;\;133_1\;\;133_2\;\;133_2\;\;188
\\126^*\;\;372\;\;406
\\133_1
\end{array}\]
\[
\begin{array}{c}
188
\\43\;\;372\;\;406
\\1\;\;89\;\;133_2\;\;133_1\;\;188\;\;188
\\43\;\;372\;\;406
\\188
\end{array}\quad
\begin{array}{c}
372
\\1\;\;89\;\;89\;\;133_2\;\;133_1\;\;188
\\43\;\;372\;\;372\;\;372\;\;406
\\1\;\;89\;\;89\;\;133_2\;\;133_1\;\;188
\\372
\end{array}\quad
\begin{array}{c}
406
\\133_2\;\;133_1\;\;188
\\126_1\;\;126_2\;\;372\;\;406
\\133_2\;\;133_1\;\;188
\\406
\end{array}\]

\section{$A_{12}$, characteristic $5$}

Here the principal blocks of $A_{12}$ and $S_{12}$ are Morita equivalent, and the non-principal block of $S_{12}$ of full defect is dual to it. We give the projectives and the factors of $S^\lambda$ for the principal bock only.

\begin{center}
\begin{tabular}{ccc}
\hline $\lambda$ &$\dim(D^\lambda)$& Factors of $S^\lambda$
\\\hline $(12)$&$1$&$1$
\\ $(9,3)$ & $153$ & $153,1$
\\ $(8,3,1)$ & $738$ & $738,153$
\\ $(7,5)$ & $144$ & $144,153$
\\ $(7,4,1)$ & $372$ & $372,144,738,153,1$
\\ $(7,3,1^2)$ & $1266$ & $1266,372,738$
\\ $(7,2,1^3)$ & $462$ & $462,1266$
\\ $(6,3,1^3)$ & $1596$ & $1596,462,1266,372$
\\ $(5,3,1^4)$ & $1506$ & $1506,1596,462$
\\ $(4^3)$ & $89$ & $89,372,1$
\\ $(4,3^2,1^2)$ & $1957$ & $1957,89,1596,372,144$
\\ $(4,3,2,1^3)$ & $573$ & $573,1957,1506,1596$ 
\\ $(3^3,2,1)$ & $11$ & $11,1957,144$
\\ $(3^2,2^2,1^2)$ & $43$ & $43,11,573,1957,89$
\\\hline
\end{tabular}
\end{center}

\[
\begin{array}{c}
1
\\ 153\;\;372
\\ 1\;\;1\;\;89\;\;144
\\153\;\;372
\\1
\end{array}\qquad \begin{array}{c}
11
\\43\;\;1957
\\11\;\;11\;\;89\;\;144\;\;573
\\43\;\;1957
\\11
\end{array}\qquad \begin{array}{c}
43
\\11\;\;89\;\;573
\\43\;\;43\;\;1957
\\11\;\;89\;\;573
\\43
\end{array}\qquad \begin{array}{c}
89
\\43\;\;372\;\;1957
\\1\;\;11\;\;89\;\;89\;\;144\;\;573\;\;1596
\\43\;\;372\;\;1957
\\89
\end{array}\]
\[\begin{array}{c}
144
\\153\;\;372\;\;1957
\\1\;\;11\;\;89\;\;144\;\;144\;\;738\;\;1596
\\153\;\;372\;\;1957
\\144
\end{array}\qquad \begin{array}{c}
153
\\1\;\;144\;\;738
\\153\;\;153\;\;372
\\1\;\;144\;\;738
\\153
\end{array}\qquad \begin{array}{c}
372
\\1\;\;89\;\;144\;\;738\;\;1266\;\;1596
\\153\;\;372\;\;372\;\;372\;\;462\;\;1957
\\1\;\;89\;\;144\;\;738\;\;1266\;\;1596
\\372
\end{array}\]
\[\begin{array}{c}
462
\\1266\;\;1596
\\372\;\;462\;\;462\;\;1506
\\1266\;\;1596
\\462
\end{array}\qquad\begin{array}{c}
573
\\43\;\;1506\;\;1957
\\11\;\;89\;\;573\;\;573\;\;1596
\\43\;\;1506\;\;1957
\\573
\end{array}\qquad\begin{array}{c}
738
\\153\;\;372
\\1\;\;144\;\;738\;\;1266
\\153\;\;372
\\738
\end{array}\qquad\begin{array}{c}
1266
\\372\;\;462
\\738\;\;1266\;\;1596
\\372\;\;462
\\1266
\end{array}\]
\[\begin{array}{c}
1506
\\573\;\;1596
\\462\;\;1506\;\;1957
\\573\;\;1596
\\1506
\end{array}\qquad\begin{array}{c}
1596
\\372\;\;462\;\;1506\;\;1957
\\89\;\;144\;\;573\;\;1266\;\;1596\;\;1596
\\372\;\;462\;\;1506\;\;1957
\\1596
\end{array}\qquad\begin{array}{c}
1957
\\11\;\;89\;\;144\;\;573\;\;1596
\\43\;\;372\;\;1506\;\;1957\;\;1957
\\11\;\;89\;\;144\;\;573\;\;1596
\\1957
\end{array}\]

\bibliography{references}

\end{document}